\documentclass[11pt,letterpaper]{article}
\usepackage{amsmath,amssymb,amsthm,mathtools,bm}
\usepackage{booktabs}
\usepackage{graphicx}
\usepackage{tikz}
\usetikzlibrary{arrows.meta,positioning}
\usepackage[numbers,sort&compress]{natbib}
\usepackage[hypertexnames=false,hidelinks]{hyperref}
\usepackage{geometry}
\usepackage{placeins}
\graphicspath{{figures/}}

\newcommand{\R}{\mathbb{R}}
\newcommand{\dd}{\,\mathrm d}
\newcommand{\bx}{\boldsymbol{x}}
\newcommand{\ba}{\boldsymbol{a}}
\newcommand{\vc}{\boldsymbol{c}}

\newcommand{\ve}{\boldsymbol{e}}
\newcommand{\vz}{\boldsymbol{z}}
\newcommand{\vphi}{\boldsymbol{\phi}}
\newcommand{\vpsi}{\boldsymbol{\psi}}

\newcommand{\vxi}{\boldsymbol{\xi}}
\newcommand{\valpha}{\boldsymbol{\alpha}}
\newcommand{\dt}{\Delta t}
\newcommand{\Phit}{\widetilde{\bm\Phi}}
\newcommand{\Dt}{\widetilde{\mathbf D}}
\newcommand{\bt}{\widetilde{\boldsymbol b}}

\newcommand{\al}{\valpha}
\newcommand{\norm}[1]{\left\lVert #1\right\rVert}

\newtheorem{proposition}{Proposition}

\newtheorem{theorem}{Theorem}

\newcommand{\GeneratedSubspaceStatement}{%
For the selected assembly configurations, recomputing the weighted raw-feature SVD on the comparison rules preserves the selected assembly/comparison ranks (1D $61/61$, 2D $215/215$, 3D $260/260$). Across these configurations, the largest principal angle between the retained subspaces is $2.260\times10^{-5}$ degrees; the discarded-block and retained--discarded cross-mass spectral ratios do not exceed $9.536\times10^{-12}$ and $3.898\times10^{-16}$, respectively.
}

\newcommand{\GeneratedConfigurationTable}{%
\begin{table}[!ht]
\centering
\caption{Configurations for the assembled parabolic tests.  Here $N_2/N_1$ are the second- and first-layer widths, and $n_a/n_c$ are the numbers of Gauss--Legendre nodes per coordinate for the assembly and comparison rules.  The time steps are $\Delta t_j=\Delta t_{\max}/2^j$, $j=0,\ldots,J$.}\vspace{0.2cm}
\label{tab:configuration}
\small
\setlength{\tabcolsep}{3pt}
\begin{tabular}{lccccccc}
\toprule
Case & $N_2/N_1$ & $\gamma$ & rank & $n_a/n_c$ & $T$ & $\Delta t_{\max}$ & $J$\\
\midrule
1D & 100/160 & 28 & 61 & 260/640 & 0.24 & 0.02 & 5 \\
2D & 240/360 & 8 & 215 & 48/64 & 0.24 & 0.04 & 5 \\
3D & 260/380 & 4 & 260 & 21/24 & 0.24 & 0.04 & 4 \\
\bottomrule
\end{tabular}
\end{table}
}

\newcommand{\GeneratedQuadratureTable}{%
\begin{table}[!ht]
\centering
\caption{Quadrature perturbations for each pair of assembly and comparison rules.  Within each row, both rules use the retained basis constructed on $q_a$.  The quantity $\widehat B_{a,c}^{\rm rel}$ is the trapezoidal approximation of the estimate in Theorem~\ref{thm:quadrature-perturbation}.}\vspace{0.2cm}
\label{tab:quadrature-comparison}
\small
\setlength{\tabcolsep}{3pt}
\begin{tabular}{lccccccc}
\toprule
Case & $n_a/n_c$ & $\varepsilon_M$ & $\varepsilon_L$ & $\varepsilon_b$ & $\varepsilon_0$ & $E_{a,c}^{\rm rel}(T)$ & $\widehat B_{a,c}^{\rm rel}(T)$\\
\midrule
1D & 260/640 & $3.005\times10^{-11}$ & $2.373\times10^{-11}$ & $6.669\times10^{-12}$ & $8.758\times10^{-12}$ & $3.847\times10^{-14}$ & $3.109\times10^{-10}$ \\
2D & 26/64 & $4.693\times10^{-1}$ & $5.116\times10^{-1}$ & $2.793\times10^{-3}$ & $3.470\times10^{-6}$ & $4.913\times10^{-6}$ & $5.350\times10^{-3}$ \\
2D & 48/64 & $3.713\times10^{-7}$ & $8.164\times10^{-7}$ & $1.002\times10^{-9}$ & $5.379\times10^{-11}$ & $8.995\times10^{-12}$ & $9.741\times10^{-9}$ \\
3D & 10/24 & $8.103$ & $5.923$ & $3.636\times10^{-2}$ & $4.469\times10^{-5}$ & $1.052\times10^{-3}$ & $2.048\times10^{-1}$ \\
3D & 21/24 & $1.140\times10^{-5}$ & $1.647\times10^{-5}$ & $6.392\times10^{-8}$ & $1.821\times10^{-10}$ & $3.007\times10^{-10}$ & $6.553\times10^{-8}$ \\
\bottomrule
\end{tabular}
\end{table}
}

\newcommand{\GeneratedComparisonConvergenceTable}{%
\begin{table}[!ht]
\centering
\caption{Relative changes in the perturbation indicators under refinement of the comparison rule; $\delta_\ell$ is defined in \eqref{eq:comparison-indicator-change}.  The retained basis and assembly rule are fixed.}\vspace{0.2cm}
\label{tab:comparison-rule-convergence}
\small
\setlength{\tabcolsep}{5pt}
\begin{tabular}{lcccccc}
\toprule
Case & $n_a$ & $n_c\to n'_c$ & $\delta_M$ & $\delta_L$ & $\delta_b$ & $\delta_0$\\
\midrule
2D & 48 & $56\to 64$ & $5.524\times10^{-3}$ & $6.832\times10^{-3}$ & $2.257\times10^{-3}$ & $1.519\times10^{-4}$ \\
2D & 48 & $64\to 72$ & $2.926\times10^{-5}$ & $4.218\times10^{-5}$ & $1.211\times10^{-5}$ & $1.246\times10^{-6}$ \\
3D & 21 & $24\to 28$ & $9.137\times10^{-3}$ & $1.023\times10^{-2}$ & $1.342\times10^{-2}$ & $1.063\times10^{-2}$ \\
3D & 21 & $28\to 32$ & $1.085\times10^{-5}$ & $1.759\times10^{-5}$ & $3.066\times10^{-5}$ & $1.906\times10^{-5}$ \\
\bottomrule
\end{tabular}
\end{table}
}

\newcommand{\GeneratedTemporalTable}{%
\begin{table}[!ht]
\centering
\caption{Temporal refinement on the fixed retained systems.  Each panel gives the complete coefficient-error sequence and consecutive EOCs.  Its heading reports the estimated semidiscrete solution error $e_{\rm sd}$, and its last row gives the relative total solution error $e_{\rm total}$ at the finest step.  The errors are defined in \eqref{eq:relative-errors}.}\vspace{0.2cm}
\label{tab:pde-temporal}
\footnotesize
\setlength{\tabcolsep}{2pt}
\begin{tabular}{c*{5}{cc}}
\toprule
\multicolumn{11}{l}{1D: $e_{\rm sd}=2.787\times10^{-4}$}\\[2pt]
$\dt$ & \multicolumn{2}{c}{implicit midpoint} & \multicolumn{2}{c}{BDF2} & \multicolumn{2}{c}{BDF3} & \multicolumn{2}{c}{Radau IIA} & \multicolumn{2}{c}{Gauss--Legendre} \\
 & $e_\alpha$ & EOC & $e_\alpha$ & EOC & $e_\alpha$ & EOC & $e_\alpha$ & EOC & $e_\alpha$ & EOC \\
\midrule
$0.02$ & $9.760\times10^{-5}$ & -- & $1.427\times10^{-6}$ & -- & $2.598\times10^{-7}$ & -- & $2.026\times10^{-7}$ & -- & $5.806\times10^{-5}$ & -- \\
$0.01$ & $5.924\times10^{-5}$ & 0.72 & $3.569\times10^{-7}$ & 2.00 & $2.628\times10^{-8}$ & 3.31 & $4.289\times10^{-8}$ & 2.24 & $2.164\times10^{-5}$ & 1.42 \\
$0.005$ & $1.677\times10^{-5}$ & 1.82 & $8.957\times10^{-8}$ & 1.99 & $2.764\times10^{-9}$ & 3.25 & $8.459\times10^{-9}$ & 2.34 & $1.090\times10^{-6}$ & 4.31 \\
$0.0025$ & $1.998\times10^{-6}$ & 3.07 & $2.246\times10^{-8}$ & 2.00 & $2.722\times10^{-10}$ & 3.34 & $1.510\times10^{-9}$ & 2.49 & $1.645\times10^{-9}$ & 9.37 \\
$0.00125$ & $4.919\times10^{-7}$ & 2.02 & $5.634\times10^{-9}$ & 2.00 & $2.778\times10^{-11}$ & 3.29 & $2.414\times10^{-10}$ & 2.64 & $1.025\times10^{-10}$ & 4.00 \\
$0.000625$ & $1.230\times10^{-7}$ & 2.00 & $1.418\times10^{-9}$ & 1.99 & $1.287\times10^{-11}$ & 1.11 & $3.514\times10^{-11}$ & 2.78 & $6.403\times10^{-12}$ & 4.00 \\
\addlinespace[2pt]
$e_{\rm total}$ & \multicolumn{2}{c}{$2.787\times10^{-4}$} & \multicolumn{2}{c}{$2.787\times10^{-4}$} & \multicolumn{2}{c}{$2.787\times10^{-4}$} & \multicolumn{2}{c}{$2.787\times10^{-4}$} & \multicolumn{2}{c}{$2.787\times10^{-4}$} \\
\midrule\addlinespace[3pt]
\multicolumn{11}{l}{2D: $e_{\rm sd}=2.274\times10^{-5}$}\\[2pt]
$\dt$ & \multicolumn{2}{c}{implicit midpoint} & \multicolumn{2}{c}{BDF2} & \multicolumn{2}{c}{BDF3} & \multicolumn{2}{c}{Radau IIA} & \multicolumn{2}{c}{Gauss--Legendre} \\
 & $e_\alpha$ & EOC & $e_\alpha$ & EOC & $e_\alpha$ & EOC & $e_\alpha$ & EOC & $e_\alpha$ & EOC \\
\midrule
$0.04$ & $3.859\times10^{-3}$ & -- & $1.320\times10^{-2}$ & -- & $3.595\times10^{-3}$ & -- & $5.457\times10^{-4}$ & -- & $6.384\times10^{-5}$ & -- \\
$0.02$ & $9.568\times10^{-4}$ & 2.01 & $3.595\times10^{-3}$ & 1.88 & $1.269\times10^{-4}$ & 4.82 & $7.172\times10^{-5}$ & 2.93 & $4.146\times10^{-6}$ & 3.94 \\
$0.01$ & $2.387\times10^{-4}$ & 2.00 & $9.021\times10^{-4}$ & 1.99 & $1.090\times10^{-5}$ & 3.54 & $9.261\times10^{-6}$ & 2.95 & $2.379\times10^{-7}$ & 4.12 \\
$0.005$ & $5.965\times10^{-5}$ & 2.00 & $2.235\times10^{-4}$ & 2.01 & $2.425\times10^{-6}$ & 2.17 & $1.179\times10^{-6}$ & 2.97 & $1.466\times10^{-8}$ & 4.02 \\
$0.0025$ & $1.491\times10^{-5}$ & 2.00 & $5.551\times10^{-5}$ & 2.01 & $3.762\times10^{-7}$ & 2.69 & $1.488\times10^{-7}$ & 2.99 & $9.146\times10^{-10}$ & 4.00 \\
$0.00125$ & $3.728\times10^{-6}$ & 2.00 & $1.382\times10^{-5}$ & 2.01 & $5.162\times10^{-8}$ & 2.87 & $1.869\times10^{-8}$ & 2.99 & $5.714\times10^{-11}$ & 4.00 \\
\addlinespace[2pt]
$e_{\rm total}$ & \multicolumn{2}{c}{$2.305\times10^{-5}$} & \multicolumn{2}{c}{$2.660\times10^{-5}$} & \multicolumn{2}{c}{$2.274\times10^{-5}$} & \multicolumn{2}{c}{$2.274\times10^{-5}$} & \multicolumn{2}{c}{$2.274\times10^{-5}$} \\
\midrule\addlinespace[3pt]
\multicolumn{11}{l}{3D: $e_{\rm sd}=5.183\times10^{-5}$}\\[2pt]
$\dt$ & \multicolumn{2}{c}{implicit midpoint} & \multicolumn{2}{c}{BDF2} & \multicolumn{2}{c}{BDF3} & \multicolumn{2}{c}{Radau IIA} & \multicolumn{2}{c}{Gauss--Legendre} \\
 & $e_\alpha$ & EOC & $e_\alpha$ & EOC & $e_\alpha$ & EOC & $e_\alpha$ & EOC & $e_\alpha$ & EOC \\
\midrule
$0.04$ & $9.038\times10^{-4}$ & -- & $1.571\times10^{-4}$ & -- & $2.115\times10^{-4}$ & -- & $2.299\times10^{-5}$ & -- & $5.719\times10^{-6}$ & -- \\
$0.02$ & $2.238\times10^{-4}$ & 2.01 & $6.650\times10^{-5}$ & 1.24 & $2.399\times10^{-5}$ & 3.14 & $3.168\times10^{-6}$ & 2.86 & $3.264\times10^{-7}$ & 4.13 \\
$0.01$ & $5.581\times10^{-5}$ & 2.00 & $1.957\times10^{-5}$ & 1.76 & $2.842\times10^{-6}$ & 3.08 & $4.148\times10^{-7}$ & 2.93 & $2.039\times10^{-8}$ & 4.00 \\
$0.005$ & $1.394\times10^{-5}$ & 2.00 & $5.235\times10^{-6}$ & 1.90 & $3.441\times10^{-7}$ & 3.05 & $5.304\times10^{-8}$ & 2.97 & $1.277\times10^{-9}$ & 4.00 \\
$0.0025$ & $3.485\times10^{-6}$ & 2.00 & $1.350\times10^{-6}$ & 1.95 & $4.228\times10^{-8}$ & 3.02 & $6.705\times10^{-9}$ & 2.98 & $7.980\times10^{-11}$ & 4.00 \\
\addlinespace[2pt]
$e_{\rm total}$ & \multicolumn{2}{c}{$5.195\times10^{-5}$} & \multicolumn{2}{c}{$5.185\times10^{-5}$} & \multicolumn{2}{c}{$5.183\times10^{-5}$} & \multicolumn{2}{c}{$5.183\times10^{-5}$} & \multicolumn{2}{c}{$5.183\times10^{-5}$} \\
\bottomrule
\end{tabular}
\end{table}
}

\providecommand{\GeneratedSubspaceStatement}{}

\begin{document}
\emergencystretch=2em
\hypersetup{pdfauthor={Xi Zhang, Lili Ju, Jin Zhao},pdftitle={Quadrature Perturbations and High-Order Time Discretization for Prescribed-Feature Evolutionary GTransNet Galerkin Systems}}
\title{Quadrature Perturbations and High-Order Time Discretization for Prescribed-Feature Evolutionary GTransNet Galerkin Systems}
\author{Xi Zhang\thanks{School of Mathematical Sciences, Chengdu University of Technology,
Chengdu, Sichuan 610059, China. Email: \href{mailto:zxi@cdut.edu.cn}{zxi@cdut.edu.cn}.}
\and
Lili Ju\thanks{Department of Mathematics, University of South Carolina,
Columbia, SC 29208, USA. Email: \href{mailto:ju@math.sc.edu}{ju@math.sc.edu}.}
\and
Jin Zhao\thanks{Corresponding author. Academy for Multidisciplinary Studies,
Capital Normal University, Beijing 100048, China. Email:
\href{mailto:zjin@cnu.edu.cn}{zjin@cnu.edu.cn}.}}
\date{}
\maketitle

\begin{abstract}
In this paper we study quadrature perturbations and high-order time discretization for a prescribed-feature evolutionary generalized transferable neural network (Evo-GTransNet) Galerkin method for linear self-adjoint parabolic equations. Assembly quadrature, together with truncated singular value decomposition, is used to select a numerically resolved trial space and obtain coordinates with an identity mass matrix. Since this normalization depends on the assembly quadrature, the retained functions are also evaluated using a separate comparison quadrature, which is finer in the numerical experiments. For the two resulting semidiscrete systems, we derive an exact equation for the difference between their trajectories and an energy estimate that separates errors arising from the effective operator, load vector, and initial projection. The estimate is conditional on the trial space selected by the assembly quadrature. The retained rank, discarded mass, cross mass, and principal angles are used to assess this selection step. Refining the comparison quadrature yields stable retained spaces for the reported configurations, whereas deliberate underintegration leads to much larger perturbations. On the fixed retained systems, the temporal convergence of implicit midpoint, BDF, Radau IIA, and Gauss--Legendre methods is consistent with their classical orders over the resolved refinement ranges. The experiments distinguish the retained-coordinate error from the total solution error; the latter is measured against the PDE solution in the $L^2$ norm, with the error integrals evaluated by sufficiently fine quadrature. They also illustrate the importance of suitable starting values for BDF3 and show how the performance of IMEX-BDF schemes depends on the operator splitting and time step size.
\end{abstract}


\section{Introduction}
We distinguish two classes of neural-network discretizations for partial differential equations (PDEs).  Training-based methods, including the deep Ritz method, the deep Galerkin method, and physics-informed neural networks, determine hidden and output parameters simultaneously through nonlinear optimization~\cite{e2018deepritz,sirignano2018dgm,raissi2019pinn,karniadakis2021piml}.  Prescribed-feature methods instead prescribe nonlinear trial functions and determine only their output coefficients from collocation, least-squares, or variational equations.  Extreme learning machines (ELMs), kernel random features, and random-feature PDE solvers belong to this second class~\cite{huang2006elm,rahimi2007random,chen2022rfm}.  For a fixed realization, the trial space is finite dimensional, and Galerkin projection gives a semidiscrete evolution system.

TransNet is a geometry-informed prescribed-feature construction for PDE approximation~\cite{zhang2023transnet}.  Its single hidden layer is formed from hyperbolic-tangent transitions across sampled partition hyperplanes whose orientations and offsets are distributed with respect to a ball containing the computational domain.  Once these geometric parameters are sampled, the hidden layer remains fixed and a stationary problem is reduced to a linear least-squares problem for the output coefficients.  Generalized TransNet (GTransNet) retains this TransNet-type geometric first layer and appends additional variance-controlled, bias-free hidden layers~\cite{cheng2026gtransnet}.  The additional layers generate the final prescribed feature set while controlling activation variance; they do not introduce problem-dependent training.  In the two-hidden-layer construction used here, the first layer has width $N_1$ and the second layer supplies $N=N_2$ final trial functions.

TransNet extensions include interface methods and weak formulations.  Multi-TransNet assigns separate TransNet trial spaces to subdomains of an elliptic interface problem; Weak-TransNet uses TransNet trial functions and local radial-basis test functions in a least-squares Petrov--Galerkin formulation~\cite{lu2025multitransnet,xu2026weaktransnet}.  Training-based SD$^{2}$NN learns coupled low- and high-frequency branches for multiscale elliptic solutions and is therefore distinct from the prescribed features considered here~\cite{li2023sd2nn}.  Redundancy and ill-conditioning in feature systems have motivated physics-informed orthogonal feature learning, randomized QR/SVD preconditioning, and rank-revealing filtering in localized random-feature spaces~\cite{jia2026orthogonal,chen2024randomized,vanbeek2025rrqr}.  In the present method, numerical-rank selection is instead defined by the quadrature-mass inner product of a boundary-conforming Galerkin trial space.

A preceding study introduced the evolutionary GTransNet (Evo-GTransNet) as a prescribed-feature Galerkin method of lines for parabolic equations~\cite{ju2026evogtransnet}.  A boundary factor produces a conforming spatial trial space, a quadrature-weighted SVD truncation selects its numerically resolved part, and only the retained output coefficients evolve in time.  That work separated rank truncation from coordinate orthonormalization, established an identity-mass semidiscrete system and its energy law, and analyzed implicit-midpoint time stepping.  The present paper uses the same fixed-space formulation to study quadrature perturbations and higher-order implicit and implicit--explicit time discretizations.

Neural discretizations of evolution equations may evolve the neural parameters or keep the spatial features fixed.  The evolutional deep neural network and neural Galerkin methods evolve neural parameters through tangent-space or residual-projection equations, and randomized sparse neural Galerkin schemes reduce the number of parameters updated at each time step~\cite{du2021ednn,bruna2024ng,berman2023sparse}.  Prescribed-feature alternatives include block time marching with local ELMs, global or blockwise space--time random-feature approximations, and local randomized discontinuous Galerkin formulations~\cite{dong2021localelm,chen2023timedependent,sun2024lrnndg}.  Tang et al. combine feature-based neural trial spaces with energy-dissipative Galerkin discretizations of gradient flows~\cite{tang2026featuregalerkin}. Zhou et al. use a fixed random-feature space together with a third-order IMEX--Runge--Kutta discretization and stagewise least-squares solves~\cite{zhou2026discretetime}.  Random-feature Galerkin coordinates have also been coupled with Kalman filtering for linear-PDE data assimilation~\cite{li2026rfkalman}. 

We consider initial-value problems in a fixed, boundary-conforming GTransNet space.
Let $V_N=\operatorname{span}\{\phi_1,\ldots,\phi_N\}$ denote this trial space. Galerkin projection of a linear parabolic equation gives the ODE system
\begin{equation}\label{eq:raw-system}
    \mathbf M\dot{\vc}(t)+\mathbf A(t)\vc(t)=\boldsymbol b(t),
\end{equation}
where $\vc(t)\in\R^N$ is the vector of output coefficients, $\mathbf M\in\R^{N\times N}$ is the mass matrix, $\mathbf A(t)\in\R^{N\times N}$ is the diffusion stiffness matrix, and $\boldsymbol b(t)\in\R^N$ is the load vector. With an assembly quadrature fixed, a quadrature-weighted SVD truncation selects a retained subspace of dimension $r\le N$. Galerkin testing in this subspace yields, in identity-mass coordinates $\valpha(t)\in\R^r$, the reduced ODE system
\begin{equation}\label{eq:main-ode}
    \dot\al(t)+\mathbf K(t)\al(t)=\bt(t),
\end{equation}
where $\mathbf K(t)=\mathbf K(t)^\top\succeq0$ is the symmetric positive semidefinite retained stiffness matrix in $\R^{r\times r}$ and $\bt(t)\in\R^r$ is the retained load vector. The Evo-GTransNet construction is described in~\cite{ju2026evogtransnet}.
Evaluating the same transformed trial functions with another quadrature rule generally yields a nonidentity mass matrix.  Classical variational-crime theory treats consistency errors caused by inexact discrete forms~\cite{strangfix2008,ciarlet2002fem}, and classical parabolic Galerkin analysis supplies the corresponding evolution framework~\cite{thomee2006,wheeler1973parabolic}.  Empirical interpolation and discrete empirical interpolation approximate nonaffine or nonlinear terms through selected evaluations, whereas empirical cubature and linear-programming empirical quadrature construct reduced integration rules~\cite{barrault2004eim,chaturantabut2010deim,hernandez2017,yano2019}.  For the transferred nonpolynomial trial functions considered here, however, polynomial exactness is unavailable.  We thus compare two quadrature-defined Galerkin systems in the same retained coordinates. The perturbation result applies to any fixed differentiable trial space evaluated with positive quadrature weights. Its particular relevance to GTransNet arises from the correlated transferred features and the rule-dependent orthonormalization of the retained basis.

With the retained space and assembly rule fixed, we compare the implicit midpoint rule, BDF2/3, coefficient-aware IMEX-BDF2/3, two-stage Radau IIA, and two-stage Gauss--Legendre methods~\cite{dahlquist1963,gear1971,ascher1995,ascher1997,pareschi2005,hairerwanner1996,butcher2016} for time integration of the reduced ODE system \eqref{eq:main-ode}.  This permits the temporal error, the accuracy of multistep starting values, and the effect of implicit--explicit splitting to be examined without changing the spatial trial space.  Variable-step extensions require additional step-ratio and discrete-energy arguments; relevant analyses include adaptive BDF2 and variable-step BDF3 schemes for diffusion equations~\cite{liao2021adaptivebdf2,liao2023variablebdf3}.  All experiments below use a fixed uniform step size within each refinement run.
We derive an exact equation and an energy estimate for the difference between trajectories assembled using two quadrature rules within the same retained GTransNet space. The estimate separates the contributions of the effective operator, load, and projected initial data. Since this estimate is conditional on the assembly-selected space, we examine subspace selection separately using retained rank, discarded mass, cross-mass, and principal angles under quadrature refinement. Numerical experiments in one, two, and three spatial dimensions distinguish the time-discretization error in the retained ODE system from the total solution error relative to the PDE solution, while also examining BDF3 starting values and noncommuting IMEX splittings. The analysis does not establish a random-feature spatial convergence rate or a bound uniform under simultaneous enrichment of the dictionary and retained rank.


The remainder of the paper is organized as follows. Section~\ref{sec:space-quadrature} reviews the construction of the retained GTransNet feature space and develops the two-rule quadrature perturbation analysis. Section~\ref{sec:time-integrators} presents the high-order implicit and IMEX time discretizations and summarizes their computational costs. Section~\ref{sec:fixed-system-analysis} describes the stability setting for the fixed system and separates retained-coordinate errors from total solution errors. Section~\ref{sec:numerical} reports the numerical experiments to verify our analysis results, and Section~\ref{sec:conclusion} provides concluding remarks.

\section{Retained GTransNet feature  space and quadrature perturbations}\label{sec:space-quadrature}

Let $\Omega\subset\R^d$ be a bounded Lipschitz domain and $T>0$ be a final time. We consider the linear diffusion problem with homogeneous Dirichlet boundary conditions
\begin{equation}\label{eq:model-problem}
\left\{
\begin{aligned}
 u_t-\nabla\cdot\bigl(a(\bx,t)\nabla u\bigr)&=f(\bx,t)
 &&\text{in }\Omega\times(0,T],\\
 u(\bx,t)&=0
 &&\text{on }\partial\Omega\times(0,T],\\
 u(\bx,0)&=u_0(\bx)
 &&\text{in }\Omega.
\end{aligned}
\right.
\end{equation}
The scalar diffusion coefficient satisfies $a(\bx,t)\ge a_0>0$ for a constant $a_0$, and $f$ and $u_0$ denote the forcing term and initial datum, respectively. We assume the spatial and temporal regularity required below. With $(\cdot,\cdot)_{L^2(\Omega)}$ denoting the usual $L^2(\Omega)$ inner product, the weak formulation is to find $u(t)\in H_0^1(\Omega)$ such that, for almost every $t\in(0,T]$,
\begin{equation}\label{eq:weak-model-problem}
 (u_t,v)_{L^2(\Omega)}+
 \int_\Omega a(\bx,t)\nabla u\cdot\nabla v\,\dd\bx
 =(f,v)_{L^2(\Omega)},
 \qquad \forall v\in H_0^1(\Omega),
\end{equation}
with $u(0)=u_0$. Following the Evo-GTransNet framework, the construction below discretizes only the spatial variable. Time remains the evolution variable in the method of lines and is not appended to the network input.

\subsection{Prescribed-Feature  TransNet geometric layer and  GTransNet trial functions}
Choose $\bx_c\in\R^d$ and $R>0$ such that $\Omega\subset B_R(\bx_c)$, where $B_R(\bx_c)$ is the ball of radius $R$ centered at $\bx_c$.  Let $N_1$ be the width of the first hidden layer and let $N_2=N$ be the number of final GTransNet features.  The first layer is
\begin{equation}
 \vpsi^{[1]}(\bx)
 =\tanh\!\left(
 \mathbf\Gamma\{\mathbf A^{[1]}(\bx-\bx_c)+R\vxi\}
 \right)\in\R^{N_1}.
 \label{eq:first-layer}
\end{equation}
Here the $i$-th row of $\mathbf A^{[1]}\in\R^{N_1\times d}$ is $\ba_i^\top$ with $\|\ba_i\|_2=1$, $\vxi=(\xi_1,\ldots,\xi_{N_1})^\top\in[-1,1]^{N_1}$ contains the offsets, and $\mathbf\Gamma=\operatorname{diag}(\gamma_1,\ldots,\gamma_{N_1})$ contains positive shape parameters.  The hyperplane
\[
 \ba_i^\top(\bx-\bx_c)+R\xi_i=0
\]
is the zero level set of the $i$-th activation: $\ba_i$ fixes its orientation, and its signed distance from $\bx_c$ in the direction $\ba_i$ is $-R\xi_i$.  Thus \eqref{eq:first-layer} is the TransNet-type geometric layer \cite{zhang2023transnet}.  In the experiments, the vectors $\ba_i$ are normalized Gaussian samples, the offsets satisfy $\xi_i\stackrel{\rm iid}{\sim}\mathcal U[-1,1]$, and $\gamma_i\equiv\gamma$.

The second hidden layer is bias-free and is defined by
\begin{equation}
 \vpsi^{[2]}(\bx)
 =\tanh\!\left(\mathbf W^{[2]}\vpsi^{[1]}(\bx)\right)
 \in\R^{N_2}.
 \label{eq:second-layer}
\end{equation}
The activation $\tanh$ is applied componentwise, and $\mathbf W^{[2]}\in\R^{N\times N_1}$ has independent entries
\begin{equation}
    (\mathbf W^{[2]})_{ji}\sim
    \mathcal N\!\left(0,\frac{\delta}{N_1}\right),
    \qquad 1\le j\le N_2,\;1\le i\le N_1,
\end{equation}
where $0<\delta\le1$ is the variance-control parameter. Additional bias-free hidden layers can be constructed in the same way. We use two hidden layers, following the numerical comparisons in~\cite{cheng2026gtransnet}, which found limited additional accuracy gains from a third hidden layer.
We write $\psi_j:=\psi_j^{[2]}$, $j=1,\ldots,N:=N_2$, for the GTransNet features. All hidden parameters are sampled once and remain fixed throughout the Galerkin evolution. The analysis below uses only that the resulting trial functions are fixed and differentiable.

Let $D\in W^{1,\infty}(\Omega)\cap H_0^1(\Omega)$ be a prescribed boundary factor.  For the unit box, we take $D(\bx)=\prod_{\ell=1}^d x_\ell(1-x_\ell)$.  Define the boundary-conforming trial functions and their span by
\begin{equation}
    \phi_j(\bx)=D(\bx)\psi_j(\bx),\qquad V_N=\operatorname{span}\{\phi_1,\ldots,\phi_N\}.
\end{equation}
Then $V_N\subset H_0^1(\Omega)$.  With $\vphi(\bx)=(\phi_1(\bx),\ldots,\phi_N(\bx))^\top$ and $\vc(t)=(c_1(t),\ldots,c_N(t))^\top$, seek the semidiscrete approximation
\begin{equation}
    u_N(\bx,t)=\sum_{j=1}^N c_j(t)\phi_j(\bx)
    =\vphi(\bx)^\top\vc(t).
\end{equation}
Substitution into \eqref{eq:weak-model-problem} and Galerkin testing with $\phi_i$, $i=1,\ldots,N$, give  the ODE system \eqref{eq:raw-system}, where
\begin{equation}
\begin{aligned}
    (\mathbf M)_{ij}&=\int_\Omega \phi_j\phi_i\dd \bx,
    &
    (\mathbf A(t))_{ij}&=\int_\Omega a(\bx,t)\nabla\phi_j\cdot\nabla\phi_i\dd \bx,\quad
    b_i(t)&=\int_\Omega f(\bx,t)\phi_i(\bx)\dd\bx.
\end{aligned}
\end{equation}

\subsection{Quadrature-weighted retained basis}
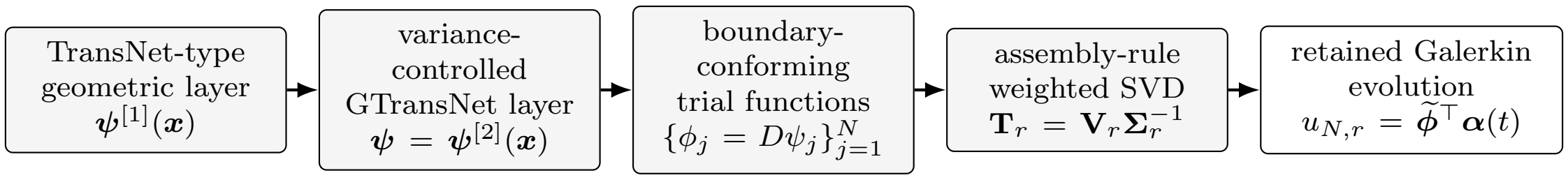
\begin{figure}[!ht]
\centering
\begin{tikzpicture}[
  fixed/.style={draw,rounded corners=2pt,fill=black!4,align=center,
                minimum height=10mm,text width=24mm,font=\scriptsize},
  evolve/.style={draw,rounded corners=2pt,align=center,
                 minimum height=10mm,text width=26mm,font=\scriptsize},
  arrow/.style={-{Latex[length=2mm]},thick},
  node distance=3mm]
\node[fixed] (first) {TransNet-type geometric layer\\$\vpsi^{[1]}(\bx)$};
\node[fixed,right=of first] (second) {variance-controlled GTransNet layer\\$\vpsi=\vpsi^{[2]}(\bx)$};
\node[fixed,right=of second] (boundary) {boundary-conforming trial functions\\$\{\phi_j=D\psi_j\}_{j=1}^N$};
\node[fixed,right=of boundary] (svd) {assembly-rule weighted SVD\\$\mathbf T_r=\mathbf V_r\bm\Sigma_r^{-1}$};
\node[evolve,right=of svd] (solution) {retained Galerkin evolution\\$u_{N,r}=\widetilde{\vphi}^{\top}\valpha(t)$};
\draw[arrow] (first) -- (second);
\draw[arrow] (second) -- (boundary);
\draw[arrow] (boundary) -- (svd);
\draw[arrow] (svd) -- (solution);
\end{tikzpicture}
\caption{Two-hidden-layer GTransNet construction and the quadrature-weighted retained basis. The hidden parameters, trial functions, and retained coordinate map are fixed; only $\valpha(t)$ evolves in time.}
\label{fig:evogtransnet-construction}
\end{figure}

Let $q_a=\{\bx_q^a,w_q^a\}_{q=1}^{Q_a}$ be the \emph{assembly} quadrature, where $Q_a$ is the total number of quadrature nodes, $\bx_q^a\in\Omega$ are the nodes and $w_q^a>0$ are the weights.  Define the feature-value matrix $\bm\Phi_a\in\R^{Q_a\times N}$ and the diagonal weight matrix $\mathbf W_a\in\R^{Q_a\times Q_a}$ by
\begin{equation}
    (\bm\Phi_a)_{qj}=\phi_j(\bx_q^a),\qquad
    \mathbf W_a=\operatorname{diag}(w_1^a,\ldots,w_{Q_a}^a).
\end{equation}
The corresponding discrete mass matrix is
\[
\mathbf M_a=\bm\Phi_a^\top\mathbf W_a\bm\Phi_a.
\]
Let $s=\operatorname{rank}(\mathbf W_a^{1/2}\bm\Phi_a)$. A compact singular value decomposition is
\begin{equation}
    \mathbf B=\mathbf W_a^{1/2}\bm\Phi_a=\mathbf U\bm\Sigma\mathbf V^\top.
\end{equation}
Here $\mathbf U\in\R^{Q_a\times s}$ and $\mathbf V\in\R^{N\times s}$ have orthonormal columns, and $\bm\Sigma=\operatorname{diag}(\sigma_1,\ldots,\sigma_s)$ with $\sigma_1\ge\cdots\ge\sigma_s>0$.  Choose a relative mass-eigenvalue threshold $\tau\in(0,1)$.  We retain the singular values satisfying $\sigma_i>\sqrt{\tau}\,\sigma_1$ and denote their number by $r\ge1$.  Let $\mathbf U_r\in\R^{Q_a\times r}$ and $\mathbf V_r\in\R^{N\times r}$ contain the corresponding singular vectors, and let $\bm\Sigma_r\in\R^{r\times r}$ contain the retained singular values.  Define
\begin{equation}\label{eq:Tr}
    \mathbf T_r=\mathbf V_r\bm\Sigma_r^{-1},\qquad
    \Phit_a=\bm\Phi_a\mathbf T_r.
\end{equation}
Thus $\mathbf T_r\in\R^{N\times r}$ maps retained coordinates to the original output coordinates, and $\Phit_a\in\R^{Q_a\times r}$ contains the retained basis evaluated on $q_a$.  This cutoff uses the same spectral-filtering mechanism as classical truncated SVD regularization; here it selects a numerically resolved trial subspace in the quadrature-mass norm rather than regularizing a single ill-posed least-squares solution~\cite{hansen1987tsvd}.  For any $\bx\in\Omega$, define the retained basis vector by
\[
    \widetilde{\vphi}(\bx)^\top=\vphi(\bx)^\top\mathbf T_r.
\]
The retained coefficient vector $\valpha(t)\in\R^r$ is related to the original output coefficient vector by $\vc(t)=\mathbf T_r\valpha(t)$, and the retained approximation is
\begin{equation}
    u_{N,r}(\bx,t)=\widetilde{\vphi}(\bx)^\top\valpha(t).
\end{equation}
Figure~\ref{fig:evogtransnet-construction} summarizes the spatial trial-space construction and retained-coefficient evolution of Evo-GTransNet~\cite{ju2026evogtransnet}.

The discarded right-singular directions are removed before time integration.  Since
\begin{equation}
    \mathbf W_a^{1/2}\Phit_a
    =\mathbf W_a^{1/2}\bm\Phi_a\mathbf V_r\bm\Sigma_r^{-1}
    =\mathbf B\mathbf V_r\bm\Sigma_r^{-1}=\mathbf U_r,
\end{equation}
we obtain
\begin{equation}\label{eq:identity-mass}
    \Phit_a^\top\mathbf W_a\Phit_a=\mathbf I_r.
\end{equation}
Here $\mathbf I_r$ denotes the $r\times r$ identity matrix. Consequently, we have
\begin{equation}
    \norm{\Phit_a\valpha}_{q_a}^{2}
    :=(\Phit_a\valpha)^\top\mathbf W_a(\Phit_a\valpha)
    =\norm{\valpha}_2^2.
\end{equation}
The retained/discarded distinction is therefore defined by the $q_a$-discrete mass norm. A discarded combination may have small mass on $q_a$ yet appreciable mass on a finer rule. The two-rule trajectory analysis below keeps $\mathbf T_r$ fixed and does not control this subspace-selection effect. We examine it numerically in Section~\ref{sec:quadrature-diagnostics}.


\subsection{Identity-mass semidiscrete evolution system}

For $1\le\ell\le d$, define the assembly-rule derivative matrices by
\begin{equation}
    (\mathbf D_{a,\ell})_{qj}=\partial_{x_\ell}\phi_j(\bx_q^a),
    \qquad
    \Dt_{a,\ell}=\mathbf D_{a,\ell}\mathbf T_r.
\end{equation}
Let
\[
 \mathbf C_a(t)=\operatorname{diag}\bigl(a(\bx_1^a,t),\ldots,a(\bx_{Q_a}^a,t)\bigr),
 \qquad
 \mathbf f_a(t)=\bigl(f(\bx_1^a,t),\ldots,f(\bx_{Q_a}^a,t)\bigr)^\top.
\]
The raw and retained assembly-rule stiffness matrices and the retained load are then respectively given by
\begin{align}
    \mathbf A_{\mathrm{raw},a}(t)
    &=\sum_{\ell=1}^d\mathbf D_{a,\ell}^\top\mathbf W_a\mathbf C_a(t)\mathbf D_{a,\ell},\\
    \mathbf K_a(t)
    &=\mathbf T_r^\top\mathbf A_{\mathrm{raw},a}(t)\mathbf T_r
      =\sum_{\ell=1}^d\Dt_{a,\ell}^\top\mathbf W_a\mathbf C_a(t)\Dt_{a,\ell},
      \label{eq:stiffness}\\
    \boldsymbol b_a(t)&=\Phit_a^\top\mathbf W_a\mathbf f_a(t).
\end{align}
Positive quadrature weights and $a(\bx_q^a,t)\ge a_0>0$ imply
$\mathbf K_a(t)=\mathbf K_a(t)^\top\succeq0$.  With
$\mathbf u_{0,a}=(u_0(\bx_1^a),\ldots,u_0(\bx_{Q_a}^a))^\top$, the assembly-rule initial projection is
\begin{equation}
    \valpha_a(0)=\Phit_a^\top\mathbf W_a\mathbf u_{0,a}.
\end{equation}


The assembly-quadrature system associated with \eqref{eq:raw-system} is
\begin{equation}
    \mathbf M_a\dot{\vc}(t)+\mathbf A_{\mathrm{raw},a}(t)\vc(t)=\boldsymbol b_{\mathrm{raw},a}(t),
\end{equation}
with $\boldsymbol b_{\mathrm{raw},a}(t)=\bm\Phi_a^\top\mathbf W_a\mathbf f_a(t)$.
Restricting $\vc(t)=\mathbf T_r\valpha(t)$ and testing only against the retained basis, by multiplication with $\mathbf T_r^\top$, gives
\begin{equation}\label{eq:retained-galerkin-system}
    \mathbf T_r^\top\mathbf M_a\mathbf T_r\,\dot\valpha(t)+\mathbf T_r^\top\mathbf A_{\mathrm{raw},a}(t)\mathbf T_r\,\valpha(t)
    =\mathbf T_r^\top\boldsymbol b_{\mathrm{raw},a}(t).
\end{equation}
Since \eqref{eq:identity-mass} gives $\mathbf T_r^\top\mathbf M_a\mathbf T_r=\mathbf I_r$, equation~\eqref{eq:retained-galerkin-system} yields \eqref{eq:main-ode} with $\mathbf K(t)=\mathbf K_a(t)$ and $\bt(t)=\boldsymbol b_a(t)$. Standard ODE integrators can therefore be applied directly.
When $\boldsymbol b_a(t)=0$, the retained semidiscrete energy satisfies
\begin{equation}\label{eq:semidiscrete-energy}
 \frac12\frac{d}{dt}\|\valpha(t)\|_2^2
 =-\valpha(t)^\top \mathbf K_a(t)\valpha(t)\le0.
\end{equation}
This also implies the decay of the assembly-rule norm of the approximate solution $u_{N,r}$.

\subsection{On a separate comparison quadrature rule}\label{sec:comparison-rule}
The identity \eqref{eq:identity-mass} is specific to $q_a$ and it does not imply normalization in the continuous $L^2$ norm or on another quadrature rule.  Let
$q_c=\{\bx_i^c,w_i^c\}_{i=1}^{Q_c}$ be a separately specified comparison quadrature rule with nodes $\bx_i^c\in\Omega$ and positive weights $w_i^c$.  Set $\mathbf W_c=\operatorname{diag}(w_1^c,\ldots,w_{Q_c}^c)$.  Without recomputing $\mathbf T_r$, evaluate the same transformed functions to form $\Phit_c\in\R^{Q_c\times r}$ and define
\begin{equation}
    (\Phit_c)_{ij}=\widetilde\phi_j(\bx_i^c),\qquad
    \mathbf G_c=\Phit_c^\top\mathbf W_c\Phit_c,\qquad
    \varepsilon_M=\norm{\mathbf G_c-\mathbf I_r}_2.
\end{equation}

\begin{proposition}[Comparison-rule norm equivalence]\label{prop:comparison-norm}
If $\varepsilon_M<1$, then every $\vz\in\R^r$ satisfies
\begin{equation}
 \sqrt{1-\varepsilon_M}\,\|\vz\|_2
 \le \|\Phit_c\vz\|_{q_c}
 \le \sqrt{1+\varepsilon_M}\,\|\vz\|_2.
\end{equation}
\end{proposition}
\begin{proof}
The symmetry of $\mathbf G_c-\mathbf I_r$ gives
$|\vz^\top(\mathbf G_c-\mathbf I_r)\vz|\le\varepsilon_M\|\vz\|_2^2$.
Since $\|\Phit_c\vz\|_{q_c}^2=\vz^\top\mathbf G_c\vz$, adding $\|\vz\|_2^2$ and taking square roots proves the result.
\end{proof}
In particular, it holds that
\begin{equation}\label{eq:comparison-spectrum}
 1-\varepsilon_M\le\lambda_{\min}(\mathbf G_c)
 \le\lambda_{\max}(\mathbf G_c)\le1+\varepsilon_M,
 \qquad
 \kappa_2(\mathbf G_c)\le\frac{1+\varepsilon_M}{1-\varepsilon_M}.
\end{equation}
Thus a small comparison-rule Gram defect controls both the norm-equivalence constants in Proposition~\ref{prop:comparison-norm} and the conditioning of $\mathbf G_c$.  When $\varepsilon_M\ge1$, Proposition~\ref{prop:comparison-norm} supplies no positive lower bound.

On $q_c$, let us define $\Dt_{c,\ell}$ by evaluating
$\partial_{x_\ell}\widetilde\phi_j$ at the comparison nodes, and let $\mathbf C_c(t)$,
$\mathbf f_c(t)$, and $\mathbf u_{0,c}$ denote the corresponding sampled coefficient, forcing,
and initial data, respectively.  The comparison-rule matrices and load are then given by
\begin{equation}
 \mathbf K_c(t)=\sum_{\ell=1}^d\Dt_{c,\ell}^\top\mathbf W_c\mathbf C_c(t)\Dt_{c,\ell},
 \qquad
 \boldsymbol b_c(t)=\Phit_c^\top\mathbf W_c\mathbf f_c(t).
\end{equation}
Assume that $\mathbf G_c$ is positive definite. The two trajectories are represented in the \emph{same retained GTransNet basis} and satisfy
\begin{align}
 \dot\valpha_a+\mathbf K_a(t)\valpha_a&=\boldsymbol b_a(t),
 &&\hspace{-1cm}\valpha_a(0)=\Phit_a^\top\mathbf W_a\mathbf u_{0,a},
 \label{eq:assembly-trajectory}\\
 \mathbf G_c\dot\valpha_c+\mathbf K_c(t)\valpha_c&=\boldsymbol b_c(t),
 &&\hspace{-1cm}\valpha_c(0)=\mathbf G_c^{-1}\Phit_c^\top\mathbf W_c\mathbf u_{0,c}.
 \label{eq:comparison-trajectory}
\end{align}
Recomputing the SVD on $q_c$ would enforce an identity $q_c$-mass matrix and generally alter the retained coordinates, and possibly the retained subspace.  It would therefore no longer isolate quadrature perturbations within the fixed space selected on $q_a$.

\subsection{Two-rule trajectory perturbation}
For $\mathbf G_c=\mathbf G_c^\top\succ0$, define $\|\vz\|_{\mathbf G_c}=(\vz^\top\mathbf G_c\vz)^{1/2}$ for $\vz\in\R^r$.
\begin{theorem}[Two-rule trajectory perturbation]\label{thm:quadrature-perturbation}
Assume that $\mathbf G_c=\mathbf G_c^\top\succ0$ and that
\[
\mathbf K_a,\mathbf K_c\in C([0,T];\R^{r\times r}),\qquad
\boldsymbol b_a,\boldsymbol b_c\in C([0,T];\R^r),
\]
with $\mathbf K_c(t)=\mathbf K_c(t)^\top\succeq0$ on $[0,T]$.
Let $\valpha_a$ and $\valpha_c$ be the corresponding continuously differentiable solutions of \eqref{eq:assembly-trajectory} and \eqref{eq:comparison-trajectory}, set $\ve=\valpha_a-\valpha_c$, and define
\begin{equation}\label{eq:quadrature-residual}
 \boldsymbol r_{a,c}(t)=\{\mathbf K_c(t)-\mathbf G_c\mathbf K_a(t)\}\valpha_a(t)
             +\mathbf G_c\boldsymbol b_a(t)-\boldsymbol b_c(t).
\end{equation}
Then, for every $t\in[0,T]$,
\begin{equation}\label{eq:quadrature-trajectory-bound}
 \|\ve(t)\|_{\mathbf G_c}
 \le \|\ve(0)\|_{\mathbf G_c}
 +\int_0^t\|\mathbf G_c^{-1/2}\boldsymbol r_{a,c}(s)\|_2\,ds
 =:B_{a,c}(t),
\end{equation}
\end{theorem}

\begin{proof}
Multiplying \eqref{eq:assembly-trajectory} by $\mathbf G_c$ and subtracting \eqref{eq:comparison-trajectory} gives 
\begin{equation}\label{eq:quadrature-error-system}
\mathbf G_c\dot{\ve}+\mathbf K_c\ve=\boldsymbol r_{a,c}.
\end{equation}
Taking the Euclidean inner product of \eqref{eq:quadrature-error-system} with $\ve$ and using $\mathbf K_c\succeq0$ gives
\[
 \frac12\frac{d}{dt}\|\ve\|_{\mathbf G_c}^2
 \le \ve^\top\boldsymbol r_{a,c}
 \le\|\ve\|_{\mathbf G_c}\|\mathbf G_c^{-1/2}\boldsymbol r_{a,c}\|_2.
\]
For $\eta>0$, set $y_\eta=(\|\ve\|_{\mathbf G_c}^2+\eta^2)^{1/2}$ and $h=\|\mathbf G_c^{-1/2}\boldsymbol r_{a,c}\|_2$.  Then
\begin{equation}\label{eq:regularized-error-bound}
 y_\eta'
 =\frac{\ve^\top\mathbf G_c\dot{\ve}}{y_\eta}
 \le \frac{\|\ve\|_{\mathbf G_c}}{y_\eta}h
 \le h,
\end{equation}
including at $\ve=\boldsymbol 0$. Integration of \eqref{eq:regularized-error-bound} on $[0,t]$ followed by $\eta\downarrow0$ proves \eqref{eq:quadrature-trajectory-bound}.
\end{proof}

The residual can be written exactly in terms of the effective comparison-rule operator and load,
\begin{equation}\label{eq:quadrature-residual-factorization}
  \boldsymbol r_{a,c}(t)=\mathbf G_c\big[(\mathbf L_c(t)-\mathbf K_a(t))\valpha_a(t)+\boldsymbol b_a(t)-\boldsymbol g_c(t)\big],
\end{equation}
with
\[
\mathbf L_c(t)=\mathbf G_c^{-1}\mathbf K_c(t),\qquad
\boldsymbol g_c(t)=\mathbf G_c^{-1}\boldsymbol b_c(t).
\]
If $\varepsilon_M<1$, Proposition~\ref{prop:comparison-norm} and \eqref{eq:quadrature-trajectory-bound}  give
\begin{align}
 \|\ve(t)\|_{\mathbf G_c}
 &\le \|\ve(0)\|_{\mathbf G_c}
 +\sqrt{1+\varepsilon_M}\int_0^t
 \big(\|\mathbf L_c(s)-\mathbf K_a(s)\|_2\,\|\valpha_a(s)\|_2
      +\|\boldsymbol b_a(s)-\boldsymbol g_c(s)\|_2\big)\,ds.\label{eq:quadrature-operator-bound}
\end{align}
Moreover, $\mathbf K_a(t)\succeq0$ implies
\begin{equation}\label{eq:assembly-state-bound}
 \|\valpha_a(t)\|_2
 \le \|\valpha_a(0)\|_2+\int_0^t\|\boldsymbol b_a(s)\|_2\,ds.
\end{equation}
This follows by testing the assembly equation with $\valpha_a$ and using the same regularization at $\valpha_a=0$ as in the proof of Theorem~\ref{thm:quadrature-perturbation}.
Equations~\eqref{eq:quadrature-operator-bound}--\eqref{eq:assembly-state-bound} connect the trajectory estimate directly to differences in the effective operator, effective load, and projected initial state.  

\section{Time integrators in retained identity-mass coordinates}\label{sec:time-integrators}
All time integration schemes in this section act on the reduced ODE system \eqref{eq:main-ode}. The spatial trial space, assembly rule, SVD truncation, and derivative matrices remain fixed. This permits the effects of time discretization, multistep initialization, and operator splitting to be examined without changing the spatial discretization. Let $t_n=n\dt$, where $\dt>0$ is the time step size, and let $\valpha^n$ be the time-discrete approximation of $\valpha(t_n)$. We write $\mathbf K^n=\mathbf K(t_n)$ and $\bt^n=\bt(t_n)$.

\subsection{Fully implicit formulas}
We begin with the classical implicit midpoint scheme used and analyzed in the Evo-GTransNet setting in~\cite{ju2026evogtransnet}. It serves as a second-order comparison method and gives
\begin{equation}\label{eq:implicit-midpoint}
 \left(\mathbf I_r+\frac\dt2\mathbf K^{n+1/2}\right)\al^{n+1}
 =\left(\mathbf I_r-\frac\dt2\mathbf K^{n+1/2}\right)\al^n+\dt\,\bt^{n+1/2},
\end{equation}
where $t_{n+1/2}=t_n+\dt/2$, $\mathbf K^{n+1/2}=\mathbf K(t_{n+1/2})$, and $\bt^{n+1/2}=\bt(t_{n+1/2})$.  
For a homogeneous scalar mode with eigenvalue $\lambda\ge0$, its amplification factor is $(1-\dt\lambda/2)/(1+\dt\lambda/2)$ and has modulus at most one.  Hence each midpoint step is contractive in the assembly norm. No commutativity between matrices at different time levels is needed for this one-step statement.

To compare one-step and multistep integration, we also consider BDF2 and BDF3. Their application to \eqref{eq:main-ode} gives
\begin{align}
 \left(\frac32 \mathbf I_r+\dt \mathbf K^{n+1}\right)\al^{n+1}
 &=2\al^n-\frac12\al^{n-1}+\dt\,\bt^{n+1},\label{eq:bdf2}\\
 \left(\frac{11}{6}\mathbf I_r+\dt \mathbf K^{n+1}\right)\al^{n+1}
 &=3\al^n-\frac32\al^{n-1}+\frac13\al^{n-2}
   +\dt\,\bt^{n+1}.\label{eq:bdf3}
\end{align}
The formulas \eqref{eq:bdf2} and \eqref{eq:bdf3} have formal orders two (BDF2) and three (BDF3), respectively.
For third-order convergence, BDF3 uses starting values $\valpha^1$ and $\valpha^2$ with errors of order $O(\dt^3)$. In the assembled PDE experiments, these values are computed by an adaptive Radau solve restricted to $[0,2\dt]$. We also compare two consecutive implicit-midpoint steps and two backward-Euler steps in Section~\ref{sec:numerical}. Over the shrinking interval $[0,2\dt]$, the former produce $O(\dt^3)$ starting errors, whereas the latter produce only $O(\dt^2)$ errors and reduce the observed global order.

\subsection{Coefficient-aware IMEX-BDF formulas}
For each scheme, the implicit--explicit partition below permits reuse of the implicit matrix factorization at a fixed step size. We study how this splitting affects the computed trajectory. Write $\mathbf K(t)=\mathbf K_0+\mathbf K_1(t)$, where $\mathbf K_0$ is independent of time. We treat $\mathbf K_0$ implicitly, evaluate $\mathbf K_1^{n+1}=\mathbf K_1(t_{n+1})$, and extrapolate only the unknown coefficient vector in $\mathbf K_1(t)\valpha(t)$. With
\[
 \valpha^{*,n+1}_{(2)}=2\valpha^n-\valpha^{n-1},
 \qquad
 \valpha^{*,n+1}_{(3)}=3\valpha^n-3\valpha^{n-1}+\valpha^{n-2},
\]
the IMEX-BDF2 and IMEX-BDF3 schemes take the form
\begin{align}
 \left(\frac32\mathbf I_r+\dt \mathbf K_0\right)\al^{n+1}
 &=2\al^n-\frac12\al^{n-1}
 +\dt\left[\bt^{n+1}-\mathbf K_1^{n+1}\al^{*,n+1}_{(2)}\right],\label{eq:imex2}\\
 \left(\frac{11}{6}\mathbf I_r+\dt \mathbf K_0\right)\al^{n+1}
 &=3\al^n-\frac32\al^{n-1}+\frac13\al^{n-2}
 +\dt\left[\bt^{n+1}-\mathbf K_1^{n+1}\al^{*,n+1}_{(3)}\right].\label{eq:imex3}
\end{align}
The formulas \eqref{eq:imex2} and \eqref{eq:imex3} have formal orders two (IMEX-BDF2) and three (IMEX-BDF3), respectively. Convergence at these orders requires the usual smoothness, stability, and starting-value assumptions.
When no separate high-order starting procedure is used, IMEX-BDF2 is initialized by the IMEX-Euler step
\begin{equation}\label{eq:imex-euler-start}
 (\mathbf I_r+\dt \mathbf K_0)\valpha^1
 =\valpha^0+\dt\{\bt^1-\mathbf K_1^1\valpha^0\}.
\end{equation}
Its $O(\dt^2)$ starting error is sufficient for second-order convergence.  IMEX-BDF3 requires $O(\dt^3)$ values at the first two levels.  The explicit component introduces a split-dependent restriction that cannot be inferred from the eigenvalues of $\mathbf K_0$ alone, especially when $\mathbf K_0$ and $\mathbf K_1(t)$ do not commute.  Classical IMEX multistep analyses likewise show that stability is a property of the coupled implicit--explicit recurrence and its partition, rather than of the implicit matrix alone~\cite{frank1997imex,hundsdorfer2007imex}.  Related high-order symmetrized IMEX-BDF analyses for parabolic operators also make the operator partition and starting corrections explicit~\cite{li2020symmetrizedimexbdf}. The numerical effect of the partition is examined in Section~\ref{sec:imex-sensitivity}.  

\subsection{Implicit Runge--Kutta methods}
For comparison with the implicit midpoint and BDF schemes, we also consider the classical two-stage Radau IIA and Gauss--Legendre methods, of orders three and four, respectively~\cite{butcher2016,hairerwanner1996}. These one-step methods provide higher-order comparisons without the multistep initialization required by BDF3. For an $s$-stage method, let $\boldsymbol Y_i\in\R^r$ be the $i$-th stage vector, and let $a_{ij}$, $b_i$, and $c_i$ be the Butcher coefficients.  The stage equations and update are
\begin{align}
 \boldsymbol Y_i&=\al^n+\dt\sum_{j=1}^s a_{ij}
 \left[-\mathbf K(t_n+c_j\dt)\boldsymbol Y_j+\bt(t_n+c_j\dt)\right],\label{eq:irk-stage}\\
 \valpha^{n+1}&=\valpha^n+\dt\sum_{i=1}^s b_i
 \left[-\mathbf K(t_n+c_i\dt)\boldsymbol Y_i+\bt(t_n+c_i\dt)\right].\label{eq:irk-update}
\end{align}
Their Butcher tableaux are respectively
\begin{equation}
{\rm Radau IIA:}\quad
\begin{array}{c|cc}
 \frac13 & \frac5{12} & -\frac1{12}\\[4pt]
 1 & \frac34 & \frac14\\[4pt] \hline
   & \frac34 & \frac14
\end{array}
\end{equation}
and 
\begin{equation}
\mbox{\rm Gauss--Legendre:}\quad
\begin{array}{c|cc}
 \frac12-\frac{\sqrt3}{6} & \frac14 & \frac14-\frac{\sqrt3}{6}\\[4pt]
 \frac12+\frac{\sqrt3}{6} & \frac14+\frac{\sqrt3}{6} & \frac14\\[4pt]\hline
   & \frac12 & \frac12
\end{array}.
\end{equation}
Each method requires the solution of one $2r\times2r$ block linear system per step, whereas a BDF step requires one $r\times r$ linear system. 
The two-stage Radau IIA method is L-stable, whereas the two-stage Gauss--Legendre method is A-stable but not L-stable. Their stability functions tend to zero and one, respectively, as the argument tends to negative infinity~\cite{hairerwanner1996,butcher2016}.

\subsection{Fixed-space construction and computational cost}
The methods above share the same retained-space construction but differ in the linear systems solved during time integration. The raw GTransNet basis and its derivatives are first evaluated on the quadrature rule $q_a$, and the 
quadrature-weighted SVD truncation defines $\mathbf T_r$. The transformed stiffness, load, and initial projection are then assembled in this fixed space, after which the selected time integrator advances the coefficient system. The same transformed functions are evaluated on the comparison quadrature rule  $q_c$ without recomputing $\mathbf T_r$. 

We count the principal dense matrix operations after the raw basis and its derivatives have been evaluated. A thin SVD of the $Q_a\times N$ weighted feature matrix has cost
\begin{equation}
 O\!\left(Q_aN\min\{Q_a,N\}\right).
\end{equation}
Transforming the basis and derivative matrices costs $O((d+1)Q_aNr)$, while direct assembly of one retained stiffness matrix costs $O(dQ_ar^2)$. Given the sampled data, assembling a retained load vector or the initial projection costs $O(Q_ar)$.

Each step of implicit-midpoint or fully implicit BDF schemes requires the solution of an $r\times r$ linear system. A dense factorization costs $O(r^3)$, followed by a solve of cost $O(r^2)$ for one right-hand side. When the step matrix varies with time, direct factorization at each step incurs the former cost repeatedly. For IMEX-BDF2 and IMEX-BDF3, the matrices to be factorized are $\frac32\mathbf I_r+\dt\mathbf K_0$ and $\frac{11}{6}\mathbf I_r+\dt\mathbf K_0$, respectively. Since $\mathbf K_0$ is independent of time, each matrix can be factorized once at a fixed step size, and its factors reused at subsequent steps. The accompanying implementation does not reuse these factorizations.

For either two-stage implicit Runge--Kutta method, a direct solution of the stage equations requires factoring a $2r\times2r$ block matrix at cost $O((2r)^3)$, followed by a solve at cost $O((2r)^2)$. This block matrix depends on the stiffness matrices at the stage times and generally changes from step to step for the time-dependent problems considered here.

In the parabolic test problems of Section~\ref{sec:numerical}, each diffusion coefficient is a sum of a constant and a product of a spatial function and a time-dependent scalar. The retained stiffness matrix is therefore a linear combination of two time-independent matrices. These components can be assembled once and combined with $O(r^2)$ work at each evaluation. This $O(r^2)$ cost is for evaluating the stiffness matrix and is additional to the linear-system solution costs above.

\section{Temporal accuracy and error decomposition}\label{sec:fixed-system-analysis}
Given a fixed spatial discretization, time-step refinement measures convergence to the retained semidiscrete solution. We first recall the accuracy and stability properties relevant to this convergence, and then relate the retained-coordinate error to the $L^2$ error in the reconstructed solution.

\subsection{Temporal accuracy and stability on the fixed retained space}
Throughout this section, the dictionary, retained rank, assembly rule $q_a$, and transformation $\mathbf T_r$ are fixed. The assembled stiffness matrix $\mathbf K(t)$ and load $\bt(t)$ may still depend on time. 
Under sufficient smoothness of $\mathbf K$, $\bt$, and $\valpha$, stability on $[0,T]$ with bounds independent of $\dt$, and starting errors of order $O(\dt^p)$, classical ODE error estimates give
\[
 \max_{0\le t_n\le T}\|\valpha(t_n)-\valpha^n\|_2\le C\dt^p
\]
for sufficiently small step sizes~\cite{hairernorsettwanner1993,hairerwanner1996,gear1971}. Here $p=2$ for implicit midpoint, BDF2, and IMEX-BDF2; $p=3$ for BDF3, IMEX-BDF3, and two-stage Radau IIA; and $p=4$ for two-stage Gauss--Legendre. The constant $C$ is independent of $\dt$, but may depend on the fixed retained system, the solution, the final time, and the IMEX partition. The estimate concerns convergence to the semidiscrete solution. Spatial approximation and quadrature errors enter when that solution is compared with the PDE solution.

The stability requirement depends on the time integrator. For BDF3, a classical modal argument applies when the retained stiffness matrix $\mathbf K$ is constant and symmetric positive definite. Let us write $\mathbf K=\mathbf Q\boldsymbol\Lambda\mathbf Q^\top$, where $\mathbf Q$ is orthogonal and $\boldsymbol\Lambda$ is the diagonal matrix of eigenvalues. The difference between two BDF3 solutions with the same discrete load satisfies the homogeneous recurrence associated with \eqref{eq:bdf3}. In the coordinates defined by $\mathbf Q$, each component satisfies a scalar BDF3 recurrence. For a component $\beta^n$ corresponding to an eigenvalue $\lambda\ge0$, substitution of $\beta^n=\zeta^n$ in \eqref{eq:bdf3} gives its characteristic polynomial
\begin{equation}\label{eq:bdf3-poly}
 \left(\frac{11}{6}+s\right)\zeta^3-3\zeta^2+\frac32\zeta-\frac13=0,
 \qquad s=\dt\lambda.
\end{equation}
The BDF3 stability region contains the negative real axis, so all roots are strictly inside the unit disk for $s>0$~\cite{gear1971,hairerwanner1996}. At $s=0$, it holds
\[
 11\zeta^3-18\zeta^2+9\zeta-2=(\zeta-1)(11\zeta^2-7\zeta+2),
\]
where $\zeta=1$ is simple and the conjugate pair has modulus $\sqrt{2/11}$. Thus perturbations in the starting values remain bounded under this autonomous recurrence and decay in modes with $\lambda>0$. This multistep stability statement does not require the solution difference to decrease in norm at every time step.

For comparison, implicit midpoint is nonexpansive at each time step in the assembly norm: the difference between two solutions with the same discrete load cannot increase in norm when the midpoint stiffness is symmetric positive semidefinite, as follows from \eqref{eq:implicit-midpoint}. This property also holds for time-dependent $\mathbf K(t)$. BDF3 is not A-stable~\cite{dahlquist1963}, and the modal argument above does not establish stability for general noncommuting $\mathbf K(t)$. For the IMEX-BDF schemes, the explicit component introduces an additional dependence on the operator partition, as discussed in Section~\ref{sec:time-integrators}. The time-dependent assembled systems and the effect of the IMEX partition are also examined numerically in Section~\ref{sec:numerical}.

\subsection{Retained-coordinate and total solution errors}\label{subsec:error-decomposition}
We now relate the temporal error in the retained coordinates to the error in the reconstructed solution. Let $u_{N,r}(\cdot,t)$ be the retained semidiscrete solution with coefficient vector $\valpha(t)$, and define the time-discrete reconstruction by
\[
 u_{N,r}^n(\bx)=\widetilde{\vphi}(\bx)^\top\valpha^n.
\]
With $\ve_n:=\valpha(t_n)-\valpha^n$, the identity-mass normalization \eqref{eq:identity-mass} gives, on the assembly rule,
\begin{equation}\label{eq:assembly-temporal-identity}
 \|\Phit_a(\valpha(t_n)-\valpha^n)\|_{q_a}
 =\|\valpha(t_n)-\valpha^n\|_2.
\end{equation}
This identity is specific to the assembly rule and does not imply normalization in the continuous $L^2(\Omega)$ norm. Since the retained basis is fixed, convergence of the coefficient vector implies convergence of the reconstructed solution in $L^2(\Omega)$.

The total solution error is measured against the PDE solution in the $L^2(\Omega)$ norm. Adding and subtracting $u_{N,r}(\cdot,t_n)$ and applying the triangle inequality gives
\begin{equation}\label{eq:solution-error-decomposition}
\begin{aligned}
 \|u(\cdot,t_n)-u_{N,r}^n\|_{L^2(\Omega)}
 &\le \underbrace{\|u(\cdot,t_n)-u_{N,r}(\cdot,t_n)\|_{L^2(\Omega)}}_{E_{\rm space/assembly}(t_n)}
 \quad+\underbrace{\|u_{N,r}(\cdot,t_n)-u_{N,r}^n\|_{L^2(\Omega)}}_{E_{\rm temporal}(t_n)}.
\end{aligned}
\end{equation}
The first term includes errors from the spatial approximation, SVD truncation, initial projection, and quadrature assembly. Note that no random-feature spatial rate is asserted. The second term is the temporal error in the reconstructed solution. The reverse triangle inequality also gives
\begin{equation}
 \left|\|u(\cdot,t_n)-u_{N,r}^n\|_{L^2(\Omega)}
       -E_{\rm space/assembly}(t_n)\right|
 \le E_{\rm temporal}(t_n).
\end{equation}
Consequently, at any fixed evaluation time, the total solution error converges to the spatial/assembly contribution as the temporal error tends to zero. This limiting value need not be a lower bound at finite step sizes, since the spatial and temporal error vectors may partially cancel. The retained-coordinate error $\|\ve_n\|_2$ measures convergence to the same semidiscrete solution and contains no spatial/assembly contribution.

\section{Numerical experiments}\label{sec:numerical}
The experiments use prescribed exact solutions to distinguish quadrature perturbations, temporal errors in the retained ODE, and total solution errors.  After specifying the numerical setup and test problems, we examine quadrature sensitivity and temporal convergence, followed by the effects of BDF3 starting values and explicitly treated diffusion.

\subsection{Experimental setup}

\subsubsection{Quadrature rules and error measures} 
The assembly and comparison rules $q_a$ and $q_c$ are tensor-product Gauss--Legendre rules. Their numbers of nodes per coordinate are denoted by $n_a$ and $n_c$, respectively. Table~\ref{tab:configuration} gives these counts for all three test problems. 
All dictionaries are generated by using $\bx_c=(1/2,\ldots,1/2)$, $R=\sqrt d/2$, and $\delta=0.7$. The SVD truncation uses the relative mass-eigenvalue threshold $\tau=10^{-11}$. One fixed dictionary realization is used for each assembled PDE case.  Within each temporal refinement sequence, the spatial discretization and all quadrature rules remain fixed while $\dt$ varies.

\GeneratedConfigurationTable

For the numerical evaluation of the error decomposition in Section~\ref{subsec:error-decomposition}, we approximate the exact semidiscrete solution at the final time $T$ by a numerical reference and report relative errors. Let $\valpha_{\rm ref}(T)$ denote the reference coefficient vector and $\valpha_{\dt}(T)$ the vector computed by the selected time integrator with step size $\dt$. The corresponding reconstructed solutions are
\[
 u_{\dt}(\bx,T)=\widetilde{\vphi}(\bx)^\top\valpha_{\dt}(T),
 \qquad
 u_{\rm ref}(\bx,T)=\widetilde{\vphi}(\bx)^\top\valpha_{\rm ref}(T).
\]
We define the following relative errors:
\begin{equation}\label{eq:relative-errors}
\begin{aligned}
 e_\alpha&=\frac{\|\valpha_{\dt}(T)-\valpha_{\rm ref}(T)\|_2}{\|\valpha_{\rm ref}(T)\|_2},\\
 e_{\rm total}&=\frac{\|u_{\dt}(\cdot,T)-u(\cdot,T)\|_{L^2(\Omega)}}{\|u(\cdot,T)\|_{L^2(\Omega)}},\\
 e_{\rm sd}&=\frac{\|u_{\rm ref}(\cdot,T)-u(\cdot,T)\|_{L^2(\Omega)}}{\|u(\cdot,T)\|_{L^2(\Omega)}}.
\end{aligned}
\end{equation}
The coefficient error $e_\alpha$ measures temporal accuracy, up to the accuracy of the reference solve. It uses a different norm and normalization from the two solution errors, so their relative magnitudes are not directly comparable. The quantities $e_{\rm total}$ and $e_{\rm sd}$ use the same $L^2(\Omega)$ norm and the same normalization. The former is the relative total solution error; the latter estimates the relative error of the assembled semidiscrete solution, including spatial approximation, SVD truncation, initial projection, and quadrature assembly. The estimate also contains the remaining error of the reference time integration. All denominators in \eqref{eq:relative-errors} are nonzero for the tests considered here.

The integrals defining the $L^2$ errors are approximated by tensor-product Gauss--Legendre quadrature with $400$, $64$, and $24$ nodes per coordinate in one, two, and three dimensions, respectively. These rules are fixed throughout the temporal refinements; a further quadrature-refinement check is reported in Section~\ref{sec:temporal-convergence}.
Reference solutions of the retained ODE are computed with the adaptive fifth-order Radau IIA method implemented in SciPy (version 1.15.3), using an analytic Jacobian and relative and absolute tolerances of $10^{-13}$ and $10^{-15}$, respectively~\cite{virtanen2020scipy}.  This reference solver differs from the fixed-step two-stage, third-order Radau IIA method used in the temporal comparisons.
In the assembled PDE experiments, BDF2 is initialized by one implicit-midpoint step, BDF3 by an adaptive Radau solve restricted to $[0,2\dt]$, IMEX-BDF2 by \eqref{eq:imex-euler-start}, and IMEX-BDF3 by the same local Radau procedure used for BDF3.  The local Radau solves for starting values use relative and absolute tolerances of $10^{-11}$ and $10^{-13}$, respectively.

\begin{samepage}
For two successive temporal refinements with step sizes $\dt_i>\dt_{i+1}$, let $e_i$ and $e_{i+1}$ denote the corresponding errors at the final time $T$.  Here $i$ labels the refinement level.  The experimental order of convergence (EOC) is defined as
\begin{equation}\label{eq:pairwise-eoc}
 {\rm EOC}_i=\frac{\log(e_i/e_{i+1})}{\log(\dt_i/\dt_{i+1})}.
\end{equation}
\end{samepage}%
Pairwise EOCs are used to identify the range in which the expected temporal order is observed. Stiff pre-asymptotic behavior and split-dependent error constants can obscure the classical rates, while reference-solution error, cancellation, or roundoff may affect the finest refinements~\cite{prothero1974,hairerwanner1996}.

\subsubsection{Test problems with prescribed exact solutions}\label{sec:test-problems}
We test some problems with prescribed exact solutions to assess the total solution error independently of the retained ODE reference.  The tests solve the model linear diffusion problem \eqref{eq:model-problem}, i.e.,
\begin{equation}\label{eq:pde-template}
    u_t-\nabla\cdot(a(\bx,t)\nabla u)=f,
\end{equation}
with the homogeneous Dirichlet boundary condition on the unit interval, square, or cube.
The forcing term is obtained from the given exact solution through
\begin{equation}\label{eq:manufactured-forcing}
    f=u_t-\nabla a\cdot\nabla u-a\Delta u.
\end{equation}
 The initial coefficients are obtained by projecting $u(\cdot,0)$ on the same assembly rule.

The one-dimensional test combines an oscillatory solution with a spatially and temporally varying diffusion coefficient.  It is used both for temporal refinement and for the starting-value comparison.  We set the exact solution 
\begin{equation}
 u(x,t)=e^{-t}\sin(6\pi x)+0.30e^{-4t}\sin(14\pi x),
\end{equation}
on the domain $\Omega=(0,1)$ and the diffusion coefficient
\begin{equation}\label{eq:1d-diffusion-test}
 a(x,t)=1+0.20\theta\left(1+\frac12\sin\frac{2\pi t}{T}\right)\sin\frac{2\pi x}{\varepsilon},
\end{equation}
with $\theta=0.6$ and $\varepsilon=0.10$ for this test; hence $a(x,t)\ge0.82$.  The two solution components have different temporal decay rates.

The two- and three-dimensional tests use smooth spatial profiles multiplied by the boundary factor $D$.  Their time-dependent amplitudes provide tests of temporal accuracy on the fixed retained systems.  On the domain $\Omega=(0,1)^2$, we let $D(x,y)=x(1-x)y(1-y)$ and set
\begin{equation}
 u(x,y,t)=D(x,y)\left[c_1(t)h_1(x,y)+c_2(t)h_2(x,y)\right],
\end{equation}
with 
\[
\begin{aligned}
h_1(x,y)&=1+0.2\cos(\pi x)\sin(\pi y), & h_2(x,y)&=\sin(\pi x)\cos(\pi y),\\
c_1(t)&=e^{-t}\{1+0.15\sin(16t)\}, & c_2(t)&=0.25e^{-2t}\cos(11.2t).
\end{aligned}
\]
The diffusion coefficient is given by
\begin{equation}
a(x,y,t)=1+0.15\left(1+0.25\sin\frac{2\pi t}{T}\right)\sin(\pi x)\sin(\pi y),
\end{equation}
so $a(x,y,t)\ge1$ on the unit square.

The three-dimensional test uses $D=x(1-x)y(1-y)z(1-z)$ and  the exact solution is defined by
\begin{equation}
 u(x,y,z,t)=D\left[e^{-t}\{1+0.1\sin(8t)\}h_1
 +0.12e^{-2t}\cos(5.6t)h_2\right],
\end{equation}
on $\Omega = (0,1)^3$ with 
\begin{equation*}
 h_1=1+0.1\cos(\pi x)\sin(\pi y)\cos(\pi z), 
 \qquad
 h_2=\sin(\pi x)\cos(\pi y)\sin(\pi z).
\end{equation*}
The diffusion  coefficient in this case  is chosen to be
\begin{equation}
 a(x,y,z,t)=1+0.1\left(1+0.25\sin\frac{2\pi t}{T}\right)
 \sin(\pi x)\sin(\pi y)\sin(\pi z).
\end{equation}
Here again $a(x,y,z,t)\ge1$ on the unit cube.  These configurations test quadrature and time integration in increasing spatial dimension.  They do not assess spatial convergence under dictionary enrichment.

These problems are used for the quadrature comparisons in Section~\ref{sec:quadrature-diagnostics} and the temporal-convergence studies in Section~\ref{sec:temporal-convergence}.  The one-dimensional problem is also used to examine the influence of starting values in Section~\ref{sec:startup-sensitivity}.

\subsection{Quadrature perturbation and retained-space diagnostics}\label{sec:quadrature-diagnostics}
We first examine the effect of quadrature rule on the operators and numerical solutions in a fixed retained space.  We then refine the comparison rule to assess these measurements and finally examine whether recomputing the SVD changes the selected space.

For each pair $(q_a,q_c)$, both rules use the transform $\mathbf T_r$ constructed on $q_a$, as required by Theorem~\ref{thm:quadrature-perturbation}.  The retained basis is fixed within each row of Table~\ref{tab:quadrature-comparison}. For rows with different assembly rules $q_a$, the retained basis is constructed separately. In addition to the mass defect $\varepsilon_M$ defined in Section~\ref{sec:comparison-rule}, we measure operator, load, and initial-projection differences.  At the sampling times $\mathcal T_c=\{jT/8\}_{j=0}^{8}$, we define the following error measurements:
\begin{equation}\label{eq:quadrature-diagnostics}
\begin{aligned}
 \varepsilon_L&=\max_{t\in\mathcal T_c}
 \frac{\|\mathbf L_c(t)-\mathbf K_a(t)\|_2}
 {\|\mathbf L_c(t)\|_2},\\
 \varepsilon_b&=\max_{t\in\mathcal T_c}
 \frac{\|\boldsymbol g_c(t)-\boldsymbol b_a(t)\|_2}
 {\|\boldsymbol g_c(t)\|_2},\\
 \varepsilon_0&=
 \frac{\|\valpha_a(0)-\valpha_c(0)\|_{\mathbf G_c}}
 {\|\valpha_c(0)\|_{\mathbf G_c}}.
\end{aligned}
\end{equation}
To assess the effect on the solution, we integrate both retained systems with adaptive Radau using relative tolerance $10^{-10}$ and absolute tolerance $10^{-12}$.  These tolerances apply to the quadrature comparison, separately from the temporal-refinement reference solves.  Let $\widehat B_{a,c}(T)$ denote the numerical approximation of the bound $B_{a,c}(T)$ in \eqref{eq:quadrature-trajectory-bound}, with its integral evaluated by a 257-point trapezoidal rule along the computed trajectory.  We report
\begin{equation}
 E_{a,c}^{\rm rel}(T)=\frac{\|\ve(T)\|_{\mathbf G_c}}{\|\valpha_c(T)\|_{\mathbf G_c}},
 \qquad
 \widehat B_{a,c}^{\rm rel}(T)=\frac{\widehat B_{a,c}(T)}{\|\valpha_c(T)\|_{\mathbf G_c}}.
\end{equation}
The reported $E_{a,c}^{\rm rel}(T)$ is a difference between numerically integrated trajectories, and $\widehat B_{a,c}^{\rm rel}(T)$ is a numerical approximation of the theoretical bound.  Neither calculation certifies the continuous-time trajectory error; differences near floating-point precision should be interpreted only as evidence that the perturbation is very small.

Table~\ref{tab:quadrature-comparison} compares the computed relative endpoint differences with $\widehat B_{a,c}^{\rm rel}(T)$. Deliberate underintegration in two and three dimensions produces substantially larger operator defects and endpoint changes than the corresponding assembly configurations in Table~\ref{tab:configuration}. In every case listed, the numerically evaluated estimate exceeds the corresponding measured difference.  Thus insufficient assembly quadrature can alter the semidiscrete trajectory even when the time integration is accurate.
The deliberately underintegrated three-dimensional case has $\varepsilon_M>1$, so Proposition~\ref{prop:comparison-norm} does not apply to that row.  The computed smallest eigenvalue is nevertheless positive, $\lambda_{\min}(\mathbf G_c)=0.5439$, with $\kappa_2(\mathbf G_c)=16.74$.  This provides numerical evidence for the positive-definiteness assumption in Theorem~\ref{thm:quadrature-perturbation}.

\GeneratedQuadratureTable

We next check whether further refinement of $q_c$ changes the measured perturbations.  This comparison keeps both $q_a$ and the resulting retained basis fixed.  For two successive node counts $n_c<n'_c$, define
\begin{equation}\label{eq:comparison-indicator-change}
 \delta_\ell(n_c,n'_c)
 =\frac{|\varepsilon_\ell(n'_c)-\varepsilon_\ell(n_c)|}
 {|\varepsilon_\ell(n'_c)|},
 \qquad \ell\in\{M,L,b,0\},
\end{equation}
where $\varepsilon_\ell(n_c)$ denotes the corresponding indicator computed with $n_c$ comparison nodes per coordinate.  All four denominators are nonzero here.  Table~\ref{tab:comparison-rule-convergence} reports these relative changes for the two- and three-dimensional tests.  The indicators stabilize on the final refinement in both dimensions.  In three dimensions, the selected comparison rule has $n_c=24$, whereas the two finest rules have $n_c=28$ and $32$.  Comparing $n_c=24$ directly with $32$ changes each indicator by at most $1.35\%$; the scale of the reported perturbations is preserved.

\GeneratedComparisonConvergenceTable

Finally, for the selected assembly configurations in Table~\ref{tab:configuration}, we examine whether recomputing the weighted SVD on $q_c$ changes the retained rank or subspace.  This check concerns space selection rather than the two trajectories compared above.  Let
$(\bm\Phi_c)_{ij}=\phi_j(\bx_i^c)$ and
$\mathbf M_c^{\rm raw}=\bm\Phi_c^\top\mathbf W_c\bm\Phi_c$.
Let $\mathbf V_\perp\in\R^{N\times(N-r)}$ have orthonormal columns spanning the Euclidean orthogonal complement of $\operatorname{Range}(\mathbf V_r)$. This complement contains the threshold-discarded right singular directions and the nullspace of the weighted feature matrix $\mathbf B=\mathbf W_a^{1/2}\bm\Phi_a$. The latter consists of coefficient vectors whose associated linear combinations of the original trial functions vanish at all assembly nodes. We measure the discarded-block and retained--discarded cross mass by
\begin{equation}\label{eq:subspace-mass-diagnostics}
 \rho_{\rm disc}
 =\frac{\|\mathbf V_\perp^\top\mathbf M_c^{\rm raw}\mathbf V_\perp\|_2}
 {\|\mathbf M_c^{\rm raw}\|_2},\qquad
 \rho_{\rm cross}
 =\frac{\|\mathbf V_r^\top\mathbf M_c^{\rm raw}\mathbf V_\perp\|_2}
 {\|\mathbf M_c^{\rm raw}\|_2}.
\end{equation}
Here $\rho_{\rm disc}$ measures the mass of the discarded directions on $q_c$, and $\rho_{\rm cross}$ measures their coupling to the retained directions, both relative to the raw mass matrix.  Both ratios are set to zero when $r=N$. Recomputing the weighted raw-feature
SVD on $q_c$ with the same relative cutoff gives a selected rank $r_c$ and
orthonormal right singular vectors $\mathbf V_{r_c,c}\in\R^{N\times r_c}$.
We compare the retained subspaces through principal angles~\cite{davis1970rotation,wedin1972svd}. These angles measure the difference between subspaces independently of the choice of orthonormal bases. When $r_c=r$, the largest principal angle is
\begin{equation}\label{eq:retained-principal-angle}
 \vartheta_{\max}
 =\arccos\!\left(\sigma_{\min}(\mathbf V_r^\top\mathbf V_{r_c,c})\right),
\end{equation}
where $\sigma_{\min}$ denotes the smallest singular value. The angles use the Euclidean inner product in raw coefficient space and are reported in degrees. The reported configurations have $r_c=r$; their rank pairs are given below.  The recomputed SVD is used only for this comparison of spaces; the trajectory calculations retain the original $\mathbf T_r$.

\GeneratedSubspaceStatement
These diagnostics indicate that the tested quadrature changes leave the retained rank unchanged and the retained subspaces close in raw coefficient space.  The one- and two-dimensional cases test the sensitivity of the truncated spaces; the three-dimensional case has $r=N$ and hence retains the full raw coefficient space.

\subsection{Temporal convergence on fixed retained systems}\label{sec:temporal-convergence}
This subsection examines convergence to the fixed semidiscrete solution as $\dt$ decreases.  Table~\ref{tab:pde-temporal} reports the coefficient error $e_\alpha$ and its EOCs, together with the estimated semidiscrete error $e_{\rm sd}$ and the finest-step total solution error $e_{\rm total}$.  The latter two quantities show how the remaining spatial and quadrature errors affect the total accuracy.

\GeneratedTemporalTable

At the two finest reported time step sizes in each test, implicit midpoint and BDF2 exhibit approximately second-order convergence, while Gauss--Legendre exhibits fourth-order convergence. The Radau IIA results in two and three dimensions and the BDF3 results in three dimensions are consistent with third-order convergence; the two-dimensional BDF3 sequence shows a trend toward this order on the final refinements. For BDF3 in one dimension, the observed order decreases on the final refinement after a range of approximately third-order convergence. This decrease occurs when the coefficient error is already of order $10^{-11}$. The one-dimensional Radau IIA EOCs increase toward three but remain below it over the reported range. These comparisons show that higher-order temporal convergence can be obtained on the fixed GTransNet trial spaces.

The retained operators depend on time.  With $\mathbf K^{(1)}=\mathbf K(T/4)$ and $\mathbf K^{(2)}=\mathbf K(3T/4)$, the normalized commutators $\|\mathbf K^{(1)}\mathbf K^{(2)}-\mathbf K^{(2)}\mathbf K^{(1)}\|_2/(\|\mathbf K^{(1)}\|_2\|\mathbf K^{(2)}\|_2)$ are $2.19\times10^{-2}$, $3.08\times10^{-3}$, and $1.72\times10^{-3}$ in one, two, and three dimensions, respectively. Thus, the retained systems cannot be decoupled into scalar equations by a time-independent change of basis. 
The final rows of Table~\ref{tab:pde-temporal} show that the finest-step total errors are close to the estimated semidiscrete errors for BDF3, Radau IIA, and Gauss--Legendre. For these three methods in all three dimensions, an additional comparison shows that $\|u_{\dt}(\cdot,T)-u_{\rm ref}(\cdot,T)\|_{L^2(\Omega)}$ is less than $0.23\%$ of $\|u_{\rm ref}(\cdot,T)-u(\cdot,T)\|_{L^2(\Omega)}$, with both norms evaluated using the same quadrature. Together with the error decomposition in Section~\ref{subsec:error-decomposition}, these measurements indicate that the estimated semidiscrete error dominates the remaining total error in these runs. The two-dimensional BDF2 total error remains visibly larger than the estimated semidiscrete error.

\begin{figure}[!ht]
\centering
\includegraphics[width=.65\textwidth]{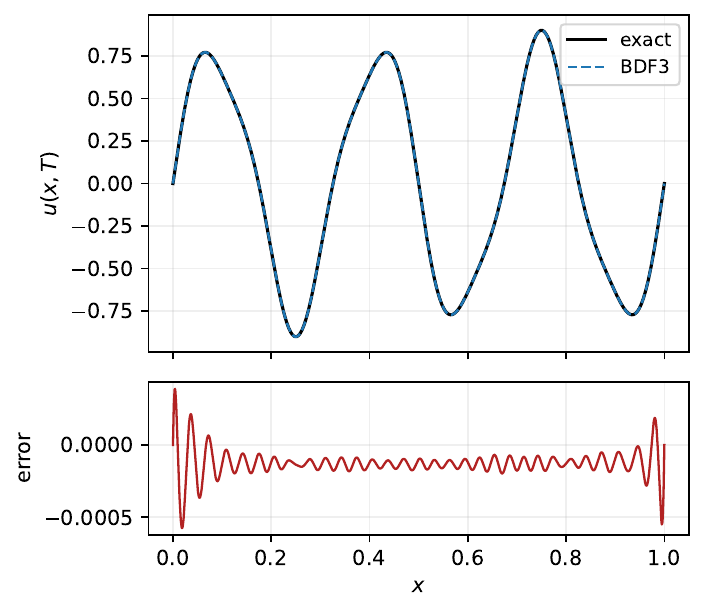}
\vspace{-0.3cm}
\caption{Exact and BDF3 profiles, with signed pointwise error, for the one-dimensional test at $T=0.24$ and $\dt=0.0025$.}
\label{fig:profile-1d}
\end{figure}

\begin{figure}[!ht]
\centering
\includegraphics[width=.99\textwidth]{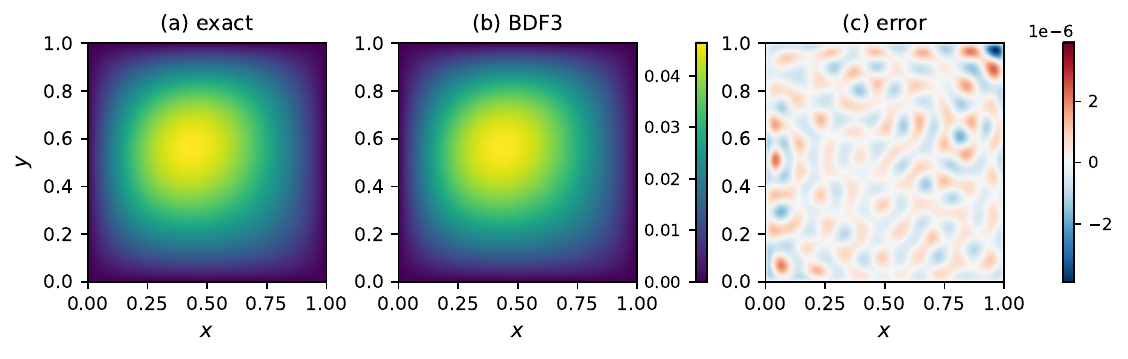}
\vspace{-0.3cm}
\caption{Exact solution, BDF3 solution, and pointwise error for the two-dimensional test at $T=0.24$ and $\dt=0.0025$.}
\label{fig:solution-snapshot}
\end{figure}

To check the numerical evaluation of the $L^2$ errors reported in Table~\ref{tab:pde-temporal}, we re-evaluate $e_{\rm sd}$ and $e_{\rm total}$ using finer quadrature rules. The number of Gauss--Legendre nodes per coordinate is increased from $400$ to $800$ in one dimension, from $64$ to $96$ in two dimensions, and from $24$ to $40$ in three dimensions. In each test, the retained basis, the reference solution at $T$, and the solutions at $T$ from all five time integrators at their finest reported time step sizes are kept fixed. All checked error values remain unchanged to the four significant digits reported.
The following plots complement the norm errors by showing the spatial distribution of the solution error.  Figure~\ref{fig:profile-1d} compares the exact and reconstructed one-dimensional profiles and their signed pointwise difference at $T=0.24$ and $\dt=0.0025$. Figure~\ref{fig:solution-snapshot} gives the corresponding two-dimensional comparison.

\subsection{Sensitivity to BDF3 starting values}\label{sec:startup-sensitivity}
\begin{figure}[!ht]
\centering
\includegraphics[width=.6\textwidth]{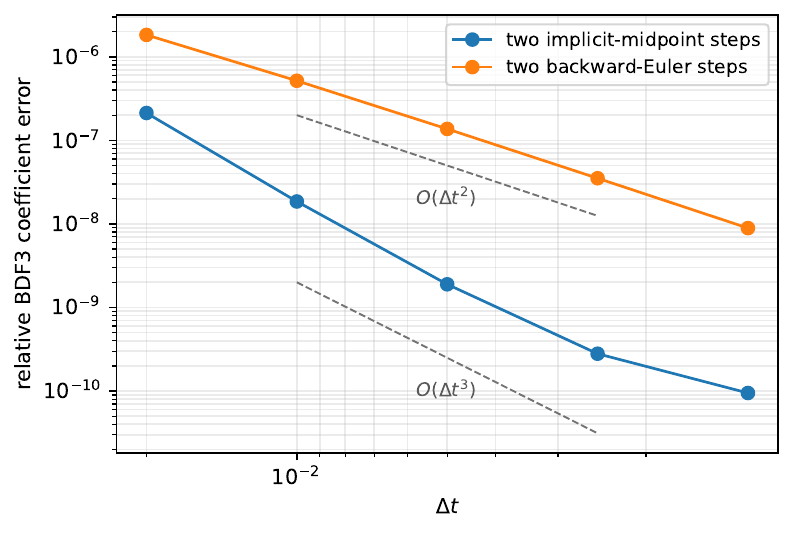}
\vspace{-0.3cm}
\caption{BDF3 coefficient errors with two implicit-midpoint or two backward-Euler starting steps, at $T=0.24$. The dashed reference lines indicate $O(\dt^3)$ and $O(\dt^2)$; their vertical positions are chosen for comparison and are not fitted to the data.}
\label{fig:startup}
\end{figure}

This experiment tests the starting-value requirement for BDF3 discussed in Section~\ref{sec:time-integrators}.  We keep the one-dimensional retained system and the subsequent BDF3 recurrence fixed, and change only the procedure used to compute the first two values.  The local adaptive-Radau baseline is already included in Table~\ref{tab:pde-temporal}; Figure~\ref{fig:startup} shows the two alternative initializations.  Every procedure uses only the interval $[0,2\dt]$.

Two consecutive implicit-midpoint steps provide $O(\dt^3)$ starting errors, which are sufficient for third-order BDF3 convergence.  The corresponding curve exhibits this rate before its finest refinement pair departs from that regime.  Two backward-Euler steps give only $O(\dt^2)$ starting errors and produce approximately second-order convergence.  The spatial system and subsequent BDF3 recurrence are identical in these computations, so this comparison isolates the effect of starting values on the observed temporal convergence.

\subsection{Effect of the explicitly treated diffusion component}\label{sec:imex-sensitivity}
This experiment examines how the strength of the explicitly treated diffusion component affects the temporal accuracy of the IMEX-BDF schemes. The main comparison concerns IMEX-BDF3, with IMEX-BDF2 included as a second-order comparison. We use a separate one-dimensional system with $N_2/N_1=100/160$, $\gamma=28$, $\delta=0.7$, $\tau=10^{-11}$, $n_a=260$, retained rank $62$, $\varepsilon=0.10$, and $T=0.20$.  The diffusion coefficient is given by \eqref{eq:1d-diffusion-test}, with $\theta\in\{0.1,0.5,1.0,1.5\}$.  The matrix $\mathbf K_0$ is assembled from the constant coefficient $a=1$, and $\mathbf K_1(t)$ from the remaining component.  Varying $\theta$ changes the diffusion coefficient and the source term in \eqref{eq:manufactured-forcing}; a separate reference trajectory is computed for each value.  The prescribed exact solution, retained space, and implicit matrix $\mathbf K_0$ remain fixed.  The experiment therefore compares problems with different explicitly treated diffusion strengths.  Both IMEX-BDF2 and IMEX-BDF3 use $\dt_j=0.02/2^j$, $j=0,\ldots,7$.  For IMEX-BDF3 with $\theta=1.5$, five additional refinements, $j=8,\ldots,12$, are included.

\begin{figure}[!htbp]
\centering
\includegraphics[width=.9\textwidth]{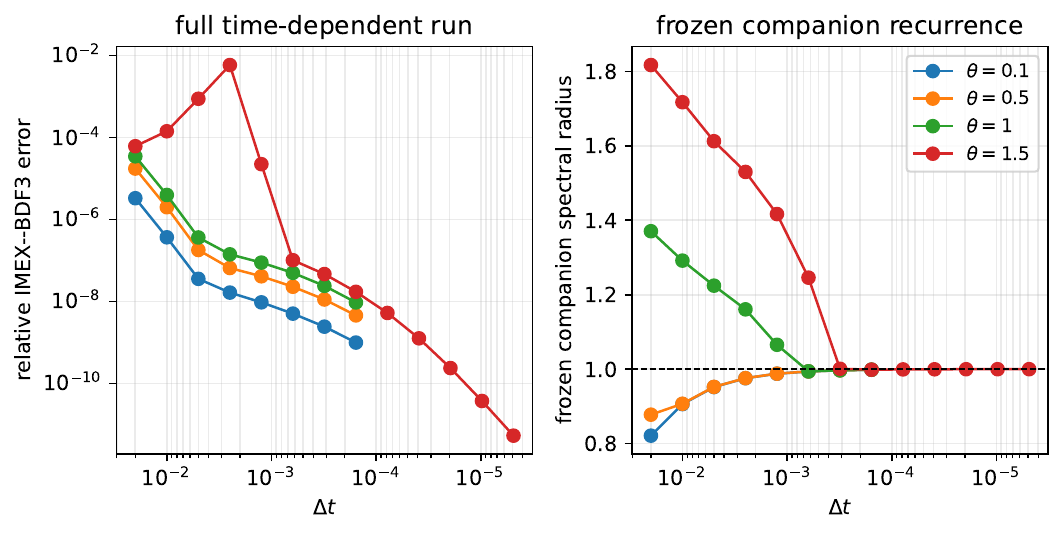}
\vspace{-0.3cm}
\caption{IMEX-BDF3 for different explicitly treated diffusion strengths $\theta$: relative coefficient error against the corresponding semidiscrete reference (left), and spectral radius of the homogeneous recurrence frozen at $t_*=T/4$ (right).}
\label{fig:imex}
\end{figure}

To examine the amplification of perturbations, we also consider the homogeneous recurrence obtained from \eqref{eq:imex3} by setting the load to zero and freezing $\mathbf K_1(t)$ at a time $t_*$. We select $t_*$ by maximizing $\|\mathbf K_1(t)\|_2$ over 101 equally spaced times in $[0,T]$; this gives $t_*=T/4$ for every tested $\theta$. For IMEX-BDF3, let $\mathbf L=\frac{11}{6}\mathbf I_r+\dt \mathbf K_0$ and form the $3r\times3r$ companion matrix whose first block row is
\begin{equation}
  \mathbf L^{-1}(3\mathbf I_r-3\dt \mathbf K_1(t_*)),\quad
  \mathbf L^{-1}(-\tfrac32\mathbf I_r+3\dt \mathbf K_1(t_*)),\quad
  \mathbf L^{-1}(\tfrac13\mathbf I_r-\dt \mathbf K_1(t_*)),
\end{equation}
with identity shift blocks below it. This matrix represents the homogeneous IMEX-BDF3 recurrence with $\mathbf K_1$ fixed at $\mathbf K_1(t_*)$. Figure~\ref{fig:imex} reports its spectral radius. For this frozen recurrence, a spectral radius below one implies decay of perturbations under repeated application of the companion matrix, whereas a radius above one permits growing modes. In the original time-dependent scheme, perturbations propagate through products of different companion matrices, so the frozen calculation does not provide a stability criterion for that scheme. 

The left panel of Figure~\ref{fig:imex} shows distinct error behavior for the tested diffusion strengths. For $\theta=0.1$, $0.5$, and $1.0$, approximately third-order convergence is observed on the first two refinements, but is not sustained on the finer grids. For $\theta=1.5$, the error initially increases as the step size decreases, reaches a maximum at $\dt=0.0025$, and then decreases. The additional refinements show a trend toward third-order convergence. 

The right panel measures the spectral radius of the frozen homogeneous recurrence.  For $\theta=1.5$, the error in the left panel starts to decrease while the frozen spectral radius in the right panel remains above one. The crossing of the unit threshold therefore does not identify the onset of error reduction in this test. The additional IMEX-BDF2 computations exhibit approximately second-order convergence on the finest refinements for all four values of $\theta$. 

\section{Concluding remarks}\label{sec:conclusion}
This work examines how quadrature assembly and time discretization affect the accuracy of prescribed-feature Evo-GTransNet Galerkin evolution.  The weighted truncated SVD selects a trial space and normalizes its mass matrix for the assembly rule; evaluation on another rule generally produces a different mass matrix.  Within the same retained space, we derive an exact trajectory-difference equation and an energy estimate that separates effective-operator, load, and initial-projection perturbations.  Retained ranks, mass-block diagnostics, and principal angles address the separate effect of quadrature on the SVD space selection.
The one-, two-, and three-dimensional experiments show that deliberate underintegration can substantially alter the assembled operators and semidiscrete trajectories, whereas the selected rules give much smaller perturbations.  On the fixed retained systems, the time integrators exhibit second-, third-, and fourth-order convergence over the respective resolved refinement ranges.  For BDF3, Radau IIA, and Gauss--Legendre at the finest reported time steps, the reconstructed temporal error is small relative to the estimated semidiscrete error.  Further time-step reduction then offers little improvement in total PDE accuracy.  The starting-value comparison isolates the influence of initialization on BDF3 convergence, while the IMEX experiments show that the observed temporal accuracy depends on both the explicitly treated diffusion component and the step size.

The resulting distinction between quadrature perturbations, temporal error, and the remaining semidiscrete error helps identify which part of the approximation requires improvement.  A higher-order time integrator can reduce the error in solving the retained ODE, but cannot by itself remove errors in the fixed spatial approximation and its assembly.  The analysis concerns fixed dictionaries and linear self-adjoint parabolic equations; convergence estimates under feature-space enrichment and extensions to adaptive or nonlinear problems remain to be developed.

\section*{Funding}
L.~Ju was supported in part by the U.S. Department of Energy under grant No.~DE-SC0025527.  J.~Zhao was supported in part by the Beijing Natural Science Foundation (No.~JR25003).

\section*{Declaration of competing interest}
The authors declare that they have no known competing financial interests or personal relationships that could have appeared to influence the work reported in this paper.

\section*{Code availability}
The source code used to generate the numerical results is provided as ancillary files with this arXiv submission.

\section*{Declaration of generative AI and AI-assisted technologies in the manuscript preparation process}
During the preparation of this work, the authors used OpenAI ChatGPT and Codex for language editing and improving readability.  The authors reviewed and edited the manuscript and take full responsibility for its content.

\bibliographystyle{plain}
\bibliography{references}

\end{document}